\documentclass[a4paper,10pt]{article}
\usepackage{amsmath,amsfonts, amssymb,verbatim}

\newcommand{\A}{\mathbb{A}}
\newcommand{\E}{\mathbb{E}}
\newcommand{\Q}{{\mathbb Q}}
\newcommand{\C}{\mathbb{C}}
\newcommand{\N}{{\mathbb N}}
\newcommand{\Z}{{\mathbb Z}}
\newcommand{\R}{{\mathbb R}}

\newcommand{\HH}{{\cal H}}

\newcommand{\F}{\mathrm{F}}
\newcommand{\U}{{\mathbb U}}
\newcommand{\V}{{\mathbb F}}

\newcommand{\X}{{\mathbb X}}
\newcommand{\Y}{\mathbb{Y}}

\newcommand{\Fp}{{\bf p}}

\newcommand{\pr}{{\rm pr}}

\newcommand{\kk}{\mathrm{k}}
\newcommand{\M}{{\bf M}}   
\newcommand{\la}{\langle}
\newcommand{\ra}{\rangle}
\newtheorem{pkt}{}[section]  
\newcommand{\bpk}{\begin{pkt}\rm }  
\newcommand{\epk}{\end{pkt}} 
\newcommand{\inv}{^{-1}}   
\newcommand{\ee}{\end{equation}}  
\newcommand{\be}{\begin{equation}}

\newcommand{\subs}{\subseteq}
\newcommand{\trd}{{\rm tr.deg}}

\newcommand{\acl}{{\rm acl}}
\newcommand{\cl}{{\rm cl}}
\newcommand{\ld}{{\rm ldim}}
\newcommand{\Mrk}{{\rm Mrk}\,}

\title{Model theory of special subvarieties  and Schanuel-type conjectures}

\author{Boris Zilber\\ University of Oxford}

\begin{document}

\maketitle

\begin{abstract}
We use the language and tools available in model theory to redefine
and clarify the rather involved notion of a special subvariety known from
the theory of Shimura varieties (mixed and pure).
\end{abstract}

\section{Introduction}
\bpk
The first part of the paper (section~\ref{s1}) is essentially a survey of  developments around the program outlined in the talk to the European Logic Colloquium-2000 and publication \cite{ZParis}. It then continues with new research which aims, on the one hand side, to extend the model-theoretic picture of \cite{ZParis} and of section~\ref{s1} to the very broad mathematical context of what we call here {\em special coverings of algebraic varieties}, and on the other hand, to use the language and the tools available in model theory to redefine and clarify the rather involved notion of a {\em special subvariety} known from the theory of Shimura varieties (mixed and pure)  and some extensions of this theory.

Our definition of special coverings of algebraic varieties includes 
semi-abelian varieties,  Shimura varieties (definitely the pure ones, and we also hope but  do not know  if the mixed ones in general satisfy all the requirement) and much more, for example, the Lie algebra covering a simple complex of Lie group
$\mathrm{SL}(2,\C).$    
    
 Recall from the discussion in \cite{ZParis} that
our specific interest in these matters arose from the connection to Hrushovski's construction of {\em new} stable structures (see e.g. \cite{Hr0}) and their relationship with  generalised Schanuel conjectures.
This subject is also closely related to
   the Trichotomy Principle and Zariski geometries. In the current paper we establish that the geometry
of an arbitrary  special covering of an algebraic variety is controlled by a Zariski geometry the closed subsets of which we call (weakly) special. The combinatorial type of  simple (i.e. strongly minimal) weakly special subsets 
 are classifiable by  the  Trichotomy Principle. Using this geometry and related dimension notions we can define a corresponding very general analogue of ``Hrushovski's predimension'' and formulate corresponding ``generalised Schanuel's conjecture'' as well as a very general forms of Andr\'e-Oort, the CIT and Pink's conjectures (Zilber-Pink conjectures). Note that in this generality one can see a considerable overlap of the generalised Schanuel  conjectures with the Andr\'e conjecture on periods \cite{An} (generalising the Grothendieck period conjecture) which prompt further questions on the model-theoretic nature of
 fundamental mathematics.

\epk

\bpk {\bf Acknowledgements.} The idea of approaching the most general setting of Schanuel-type conjectures crystallised during my 4-weeks visit at IHES at Bures-Sur-Yvette in February 2012,
and considerable
amount of work on this project was carried out during my 2-months stay at MPIM, Bonn in April-May 2012. My sincere thanks go to these two great mathematical research centres and 
people with whom 
I met and discussed the subject there. Equally crucial impact on this work comes from many conversations that I had with Jonathan Pila whose work on diophantine problems around Shimura varieties
using model-theoretic tools introduced me to the subject. Moreover, the success of his work and of the Pila-Wilkie method brought 
together people working in different areas and made it possible for me 
to consult E.Ullmo, A.Yafaev, B.Klinger and others on the topic. On model-theoretic side I used Kobi Peterzil's advice  on o-minimality issues. My thanks go to all these 
people.  
\epk
\section{Analytic and pseudo-analytic structures}\label{s1}
Recall that a strongly minimal structure $\M$ (or its theory) can be given a coarse classification by the type of {\em the combinatorial geometry} that is induced by the pregeometry $(M,\acl)$ on the set $[M\setminus \acl(\emptyset)]/_\sim,$ where $x\sim y$ iff $\acl(x)=\acl(y).$

The Trichotomy conjecture by the author stated that for any strongly minimal structure $\M$ the geometry of $\M$ is either trivial or  linear  (the two united under the name {\em locally modular}), or the geometry of $\M$ is the same as of an algebraically closed field, and in this case  the structure $\M$ is bi-interpretable with the structure of the field.

E.Hrushovski refuted this conjecture in the general setting  \cite{Hr0}.  Nevertheless  the conjecture was confirmed, by Hrushovski and the author in \cite{HZ},  for an important subclass of structures,  Zariski geometries, except for the clause stating the   bi-interpretability, where the situation turned out to be more delicate.

\bpk \label{Hr}  Recall that the main suggestion of \cite{ZParis} was to treat an (amended version 
of) Hrushovski's counter-examples as {\em pseudo-analytic } structures, analogues of classical analytic structures. Hrushovski's predimension, and the corresponding inequality $\delta(X)\ge 0,$ a key ingredient in the construction, can be seen then to directly correspond to certain type of conjectures of transcendental number theory, which we called {\em Generalised Schanuel} conjectures.

The ultimate goal in classifying the above mentioned pseudo-analytic structures has been to give a (non-first-order) axiomatisation and prove a categoricity theorem for the axiomatisable class.      
\epk
\bpk
The algebraically closed fields with pseudo-exponentiation, $\F_{\exp}=(\F;,+,\cdot, \exp),$  analogues of the classical structure   $\C_{\exp}=(\C;,+,\cdot, \exp),$ 
was the first class studied in detail. 

The axioms for $\F_{\exp}$ are as follows.

\medskip

ACF$_0:$  $\F$ is an algebraically closed  fields of characteristic  0;

\medskip

EXP: \ \ \ \ \ \  $\exp: \mathbb{G}_a(\F)\to \mathbb{G}_m(\F)$\\
is a surjective homomorphism from the additive group $\mathbb{G}_a(\F)$ to the multiplicative group $\mathbb{G}_m(\F)$ of the field
$\F$, and 
  $$\ker \exp=\omega \Z,\mbox{ for some }\omega\in \F;$$

SCH: \ \ for any finite $X,$ $$\delta(X):=\trd_\Q(X\cup \exp(X))-\ld_\Q(X)\ge 0.$$
Here $\trd(X)$ and $\ld_\Q(X)$ are the transcendence degree and the dimension of the $\Q$-linear space spanned by $X,$ dimensions of the classical pregeometries associated with the field $\F,$ and
$\delta(X)$ takes the role of {\em Hrushovski's predimension}, which gives rise to a new dimension notion and new pregeometry following Hrushovski's recipe. The inequality can be recognised as the Schanuel conjecture, if one also assumes $\F_{\exp} = \C_{\exp}.$

Since $\F_{\exp}$   results from a Hrushovski-Fr\"aisse amalgamation, the structure is 
 {\bf Existentially Closed} with respect to the embeddings respecting the predimension. This takes the form of the following property. 

{\rm EC}: For any  {\em rotund}  and {\em free} system of polynomial
equations 
$$ P(x_1,\dots,x_n,y_1,\dots,y_n)=0$$

there exists a (generic) solution satisfying $$ y_i=\exp x_i\ \ i=1,\dots,
n.$$

The term {\em rotund} has been coined by J.Kirby, \cite{Ki1}, rotund  and free  is the same as {\em normal} and free in \cite{Zexp}. We refer the reader to \cite{Ki1} or 
\cite{Zexp} for these technical definitions.

Finally we have the following
{\bf Countable Closure} property.

{\rm CC}:\ For maximal rotund systems of equations the set of solutions  is at most countable. \epk
The main result of \cite{Zexp} is the following.  
\bpk \label{main1} {\bf Theorem.} {\em Given an uncountable cardinal $\lambda,$ there is
a unique model of axioms {\rm ACF$_0$ + EXP + SCH + EC + CC} of cardinality $\lambda.$}

\epk
\bpk Recall that in \cite{ZParis} we stated: 

{\bf Conjecture.} {\em $\C_{exp}$ is isomorphic to the unique model  $\F_{\rm exp}$ of axioms \\
{\rm ACF$_0$ + EXP + SCH + EC + CC}, of cardinality continuum. }
\medskip

P. D'Aquino, A.Macintyre and G.Terzo embarked on a programme of comparing  $\F_{\rm exp}$
with  $\C_{exp},$ in particular checking the validity of (EC) which they renamed the {\em Nullstellensatz for exponential equation}.  The paper \cite{AMT1} proves that  $\F_{\rm exp}$ contains a solution 
to any equation $f(z)=0,$ where $f(z)$ is a one-variable term in $\exp,$ $+$  and $\times$ with parameters, provided
$f(z)$ is not of the form $\exp(g(z))$ for a term $g(z).$ This was proved to hold in $\C_{\rm exp}$ by W. Henson and L. Rubin using Nevanlinna Theory, \cite{HR}.   A. Schkop in \cite{Schkop} arrives at the same statement for   $\F_{\rm exp},$ deriving it directly from the axioms. 


In \cite{AMT0} P. D'Aquino, A.Macintyre and G.Terzo find a context in which the Identity Theorem of complex analysis can be checked and confirmed in $\F_{\rm exp}.$ The Identity Theorem  says in particular that if the zero set of an entire function $f(z)$
has an accumulation point then $f\equiv 0.$ Since   $\F_{\rm exp}$ comes with no topology there is no direct counterpart to the notion of accumulation in it. Instead  \cite{AMT0} considers specific subsets  of   $\F_{\rm exp},$ such that $\Q$ or the torsion points where an accumulation point can be defined in the same way as in $\C.$ In such cases, for a certain type of functions $f$     D'Aquino, Macintyre and Terzo prove that the conclusion of the Identity Theorem holds in  $\F_{\rm exp}.$ It is interesting that 
the results are obtained by invoking deep theorems of Diophantine Geometry.   

In the same spirit the present author tried, without much success, to solve the following: 
\medskip

{\bf Problem.} Describe syntactically the class of quantifier-free formulae $\sigma(z_1,\ldots,z_n)$ in the language $\{ \exp, +, \times\}$ with parameters in $\C$
for which
the subset\\ $S= \{ \bar a\in \C^n: \sigma(\bar{a})\}$  is  analytic in $\C.$

\medskip

Note that  the class is bigger than just the class of positive quantifier-free formulae. E,g, the following 
formula defines an entire function on $\C$ (and so its graph is analytic):
$$y=\left\lbrace\begin{array}{ll}
\frac{\exp(x)-1}{x}, \mbox{ if } x\neq 0,\\
e,\ \mbox{ \ \ \ \ \ \ otherwise}.
\end{array}\right.$$

\epk
\bpk 

This is a consequence of Theorems A and B.\medskip 

{\bf Theorem A.} {\em The $L_{\omega_1,\omega}(Q)$-sentence
 \begin{center}{\rm
ACF$_0$ + EXP + SCH + EC + CC}\end{center}
 is  axiomatising a {\bf quasiminimal excellent abstract elementary class} (AEC)}.
 
 \medskip

{\bf Theorem B.}  {\em A quasiminimal excellent AEC has a unique model in any
uncountable cardinality}.   \medskip

\epk
\bpk The original definition of quasiminimal excellence and the proof of Theorem B is in \cite{Zcat} based on earlier definitions and techniques  of S.Shelah.

\smallskip

{\bf Definition.} Let $\M$ be a structure for a countable language, equipped
with a pregeometry $\cl$. We say that the pregeometry of $\M$
is  {\bf quasiminimal}  if the following hold:

1. The pregeometry is determined by the language. That is, if
$\mathrm{tp}(a,\bar b) =\mathrm{tp}(a',\bar b')$ then
 $a \in \cl(\bar b)$ if and only if $a' \in \cl(\bar b').$

2. $\M$ is infinite-dimensional with respect to $\cl.$

3. (Countable closure property) If $X \subs M$ is finite, then $\cl(X)$ is countable.

4. (Uniqueness of the generic type) Suppose that $H, H' \subs M$ are countable
closed subsets, enumerated such that $\mathrm{tp}(H)=\mathrm{tp}(H').$
 If $a \in M \setminus H$ and  $a' \in M \setminus H'$ then  $\mathrm{tp}(a,H)=\mathrm{tp}(a',H').$

5. ($\aleph_0$-homogeneity over closed sets (submodels) and the empty set)
Let $H, H' \subs M$ be countable closed subsets or empty, enumerated such
that  $\mathrm{tp}(H)=\mathrm{tp}(H'),$ ),  let $\bar b, \bar b'$
be finite tuples from $M$ such that
 $\mathrm{tp}(\bar b, H)=\mathrm{tp}(\bar b', H'),$
 and let $a \in \cl(H,\bar b).$ Then there is $a' \in  M$ such
that   $\mathrm{tp}(a,\bar b, H)=\mathrm{tp}(a',\bar b', H').$

{\em Excellence} of a quasiminimal pregeometry is an extra condition on the amalgams of {\em independent systems of submodels} which we do not reproduce here in the general form but going to illustrate it in examples below.

\medskip

The proof of Theorem A relies on essential algebraic and diophantine-geometric facts and techniques and goes through an intermediate stage which is the following.

\medskip

{\bf Theorem A$_0.$} {\em  The natural  $L_{\omega_1,\omega}$-axiomatisation of the two-sorted structure $(\mathbb{G}_a(\C),\exp, \C_{\mathrm{field}})$ 
 defines a  quasiminimal excellent AEC.}

Here  $\C_{\mathrm{field}}$ is the field of complex numbers, and $\exp$ is the classical homomorphism
\be \label{hom1} \exp:\mathbb{G}_a(\C)\to  \mathbb{G}_m(\C)\ee
of the additive group onto the multiplicative group  of complex numbers.

The natural axiomatisation consists of the first-order part which consists of the theory of   $\C_{\mathrm{field}},$ the theory of  $\mathbb{G}_a(\C)$ and  the
statement that  $\exp$ is a surjective homomorphism. The only proper $L_{\omega_1,\omega}$-sentence which is added to this states that the kernel of
$\exp$ is a cyclic group. 
\medskip

These proofs were not without errors.  The original paper \cite{Zcov} of A$_0$ 
established quasiminimality of  $(\mathbb{G}_a(\C),\exp, \C_{\mathrm{field}})$ but has an error in the part proving excellence. This 
was corrected
  in   \cite{BZ}, where also a generalisation of this theorem to a positive characteristic analogue was given. But after \cite{BZ} a related error in the proof of the main Theorem A still required an extra argument which did arrive but from an unexpected direction. 
  M.Bays and J.Kirby in \cite{BK} (using \cite{Ki}), followed by  M.Bays, B.Hart, T.Hyttinen, M.Kesaala and J.Kirby \cite{BHHKK}, and further on followed by L.Haykazyan \cite{Hay} found an essential  strengthening of Theorem B, which made certain algebraic steps in 
 the proofs of Theorems A$_0$ and A redundant.  
\epk
The final result, see  \cite{BHHKK}, is the following.
\bpk {\bf Theorem B$^*.$}  {\em A quasi-minimal pregeometry can be axiomatised by an $L_{\omega_1,\omega}(Q)$-sentence which  determines an uncountably categorical class. In particular, 
the class is excellent}.   

\medskip

With the proof of this theorem the proof of the main theorem~\ref{main1} has been completed. 

Theorem B$^*$ by itself is a significant contribution to the model theory of abstract elementary classes which, remarkably, has been found while working on applications.

The significance of this model-theoretic theorem will be further emphasised in the discussion of its implications for diophantine geometry below. 
\epk
\bpk \label{d0}  Analogues of Theorem A$_0$ are now established for elliptic curves, M.Bays, and Abelian varieties by M.Bays, B.Hart and A.Pillay, based on earlier contributions by 
M.Gavrilovich, M.Bays and the author.

\medskip

{\bf Theorem A$_\mathrm{Ell}$} (M.Bays, \cite{BaPhd}). {\em Let $\E$ be an elliptic curve  without complex multiplication over a number field $\kk_0\subset \C.$     Then the natural  $L_{\omega_1,\omega}$-axiomatisation of the two-sorted structure $(\mathbb{G}_a(\C),\exp, \E(\C))$ 
 defines an uncountably categorical   AEC. In particular, this AEC is excellent.}

Here $\exp:\mathbb{G}_a(\C)\to  \E(\C)$ is the homomorphism onto the group on the elliptic curve given by the Weierstrass function and its derivative. 
The {\em natural axiomatisation} includes an axioms $\mathrm{Weil}_m$ which fixes the polynomial relation between two torsion points of order $m$ for a certain choice of $m.$ If this (first-order) axiom is dropped, the categoricity fails but not gravely -- the resulting $L_{\omega_1,\omega}$-sentence still has only finitely many (fixed number of) models in each uncountable cardinality.

The following is an extension of the previous result to
    abelian varieties incorporated in \cite{BaHaPi} by M. Bays, B. Hart and A.
    Pillay   
\medskip

{\bf Theorem A$_\mathrm{AbV}.$}  {\em Let $\A$ be an Abelian variety over a number field $\kk_0\subset \C$ and such that every endomorphism $\theta\in \mathcal{O}$ (complex multiplication) is defined over $\kk_0.$    Then the natural  $L_{\omega_1,\omega}$-axiomatisation of the two-sorted structure $(\C^g_{\mathcal{O}\cdot\mathrm{mod}},\exp, \A(\C))$ along with the first-order type of 
the kernel of $\exp$ in the two-sorted language
 defines an uncountably categorical   AEC. In particular, this AEC is excellent.}

Here $\C^g_{\mathcal{O}\cdot\mathrm{mod}}$ is the structure of the $\mathcal{O}$-module on the covering space $\C^g,$ where $g=\dim \A.$

\medskip

We explain the main ingredients of the proof. 
\epk
\bpk {\bf Definition.}  Let $\F$ be an algebraically closed field of countable transcendence degree and $B$ a finite (possibly empty) subset of its transcendence basis. 
An independent system of algebraically closed fields is
the collection $\mathcal{L} = \{ L_s\subs \F: s\subs B\}$  of algebraically closed subfields, $L_s=\acl(s).$ 

The {\bf boundary} (or the {\em crown}) ${\partial}\mathcal{L}$ of $\mathcal{L}$ is the field generated by the $L_s,$ $s\neq B,$
$${\partial}\mathcal{L}=\la L_s: s\subsetneq B\ra.$$
 
Now let us also consider a semi-Abelian variety $\A$ over a number field $\kk_0$ and assume that $\F$ is of characteristic $0,$ $\kk_0\subset \F.$

For any subfield $\kk_0\subs \kk\subs \F,$ the set $\A(\kk)$ of $\kk$-rational points of $\A$ is well-defined. And conversely, for $D\subs \A(\F)$ we write $\kk_0(D)$ the extension of $\kk_0$ by (canonical) coordinates of points of $D.$

The $\A$-boundary ${\partial}_{\A}\mathcal{L}$ of the system  $\mathcal{L}$
is the complex multiplication  submodule
$$ {\partial}_{\A}\mathcal{L}=\la \A(L_s): s\subsetneq B\ra+\mathrm{Tors}(\A),$$
where $\mathrm{Tors}(\A)$ is the torsion subgroup of $\A.$ 

The extension $\kk_{\infty}=\kk_0(\mathrm{Tors}(\A))$ will be of importance below.

We also need to use the Tate module $T(\A)$ of $\A$  which can be defined as the limit 
$$T(\A)=\lim_{\leftarrow} \A_m$$
of torsion subgroup $\A_m$ of $\A$ of order $m.$

Given $\sigma\in \mathrm{Gal}(\F:\kk_\infty)$ and $a\in \A(\kk_{\infty})$ define, for $n\in \N,$
$$\la \sigma, a\ra_n=\sigma b - b\in \A_n,$$
for an arbitrary $b\in \A$ such that $nb=a.$ This does not depend on the choice of $b$ and taking limits we have the map
$$\la\ \cdot,\,  \cdot \ra_{\infty}:\ \mathrm{Gal}(\F:\kk_\infty)\times \A(\kk_{\infty})\to T(\A),\ \ 
\la \sigma, a\ra_{\infty}=\lim_{\leftarrow}\la \sigma, a\ra_n.
$$
This also works then for $\kk\supseteq \kk_{\infty}$ in place of $ \kk_{\infty}$ and we can consider, given $a\in \A(\kk),$ the submodule 
$ \la \mathrm{Gal}(\F:\kk), a\ra_{\infty}$ of the Tate module.

The main ingredient of the proof of Theorems A$_0$ and  A$_\mathrm{AbV}$ is the following.

\medskip

{\bf Theorem}  (``Thumbtack Lemma''). {\em Let $D$ be an $\A$-boundary, let $\kk$ be a
finitely generated extension of $\kk_0(D)$ and let
$a\in \A(\kk)$ be such that  $\mathcal{O}\cdot\gamma\cap D=\{ 0\}.$
 Then $ \la \mathrm{Gal}(\F:\kk), a\ra_{\infty}$ is of finite index in $T(\A).$
 }

This splits naturally into 3 cases depending on the size $n:=|B|,$ namely cases $n=0,$ $n=1$ and $n>1.$ 

The case
 $n=0$ essentially follows from a combination of a finiteness theorem by Faltings and the Bashmakov-Ribet Kummer theory for Abelian varieties just recently finalised by M.Larsen \cite{Lar}. The similar result needed for the proof of Theorem A$_0$ for  $\mathbb{G}_m$  is just Dedekind's theory of ideals and the classical Kummer theory.  
 
 The case $n=1$ is the field of functions case of the Mordell-Weil theorem due to Lang and N\'eron, see \cite{Lan}, Theorem 6:2.
 
 The case $n>1$ is new and required a special treatment.  For $\mathbb{G}_m$ it was done in  \cite{BZ}  using the theory of specialisations (places) of fields, but  \cite{BaHaPi} finds a more direct model-theoretic argument  which covers 
also the case of Abelian varieties. 
 Using the new model-theoretic result above one has now the following, for all reasonable forms of the Thumbtack Lemma (see, in particular, \ref{Aj}).

\medskip

{\bf Corollary to Theorem B$^*$}. {\em The case  $n>1$ in the  Thumbtack Lemma
follows from the cases $n=0,1.$}
 
\epk 
\bpk \label{d1} Recall that the statement  of Theorem  A$_\mathrm{AbV}$ is weaker than that of Theorems  A$_0$ and A$_\mathrm{Ell}$ in including the complete description of the kernel of $\exp$ in the two-sorted language. The stronger version requires an extension of the   Thumbtack Lemma which includes the case $n=0\ \& \ \gamma\in \mathrm{Tors}.$ In other words one needs to characterise the action of  $\mathrm{Gal}(\tilde{\kk}:\kk)$ on $\mathrm{Tors}.$ For $\mathbb{G}_m$ this is given by the theory of cyclotomic extensions and for elliptic curves without complex multiplication by Serre's theorem.  For the general Abelian variety this is an open problem.

On the other hand  Bays, Hart and Pillay in \cite{BaHaPi} prove a broad generalisation
of  A-style theorems at the cost of fixing more parameters in the ``natural axiomatisation''. Their axioms include   the complete diagram of the prime model.
  In this setting Theorem  A holds for an arbitrary commutative
algebraic group over algebraically closed field of arbitrary characteristic. In fact,  \cite{BaHaPi} building up on \cite{BGH} by Bays, Gavrilovich and Hils
shows how to   generalise the statement to an arbitrary commutative  finite Morley rank group and proves it in this formulation.  

In terms of the Thumbtack Lemma the latter requires only $n=1$ case. This is essentially Kummer theory over  function fields in its most general, in fact model-theoretic, form.

Bays, Gavrilovich and Hils in \cite{BGH} use  this technique for an  application in algebraic geometry.
\epk

\bpk \label{d2}
Note a version of Theorem A$_\mathrm{AbV}$ formulated in  \cite{BaHaPi}.
 
 {\bf Theorem.} {\em Models of the natural first order axiomatisation of $(\C^g_{\mathcal{O}\cdot\mathrm{mod}},\exp,\A)$ are determined up to isomorphism by the transcendence degree of the field and the isomorphism type of the kernel.}
 
 This is a model-theoretic  ``decomposition'' statement for the rather complex algebraic structure, similar  to the Ax-Kochen-Ershov  ``decomposition'' of henselian valued fields into residue field and value group. 

\medskip

 
\epk
\bpk The study of the above pseudo-analytic structures shed some light on classical transcendental functions, namely the complex exponentiation $\exp$ (Theorem A, \ref{main1}), the Weierstrass  function $\mathfrak{P}(\tau_0, z)$   as a function of $z$ (Theorem A$_\mathrm{Ell}$) and more generally Abelian integrals and the corresponding exponentiation 
(Theorem A$_\mathrm{AbV}$).  Although a lot of questions remain still open, especially for the latter, a natural continuation of this program leads to questions on model theory of    $\mathfrak{P}(\tau, z)$   as a function of two variables $z,$ and similar multy-variable maps related to Abelian varieties. 

But before anything could be said about  $\mathfrak{P}(\tau, z)$ the classical function
$j(\tau),$ the modular invariant of elliptic curves $\E_{\tau},$ which can be defined in terms of  $\mathfrak{P}(\tau, z).$
\epk
\bpk \label{jaxiom} The two-sorted setting for $j(\tau)$ analogous to settings in \ref{d0} is the structure $(\HH, j, \C_{\rm field}),$ where $\HH$ is the upper half-plane as a $\mathrm{GL}^+(2,\Q)$-set, that is with the action by individual elements 
$\left(\begin{array}{ll}a\ b\\ c\ d\end{array}\right)$ from the group (rational matrices with positive determinants)
$$\left(\begin{array}{ll}a\ b\\ c\ d\end{array}\right): \tau\mapsto \frac{a\tau + b}{c\tau +d}.$$

 In the language we have names
for fixed points $t_g\in \HH$ of transformations $g,$ which are exactly the quadratic points on the upper half-plane, and the list of statements  $g t_g=t_g$ is part of the axiomatisation. The images $j(t_g)\in \C$ are called {\em special points}. They are algebraic and their values are given by the axioms of the structure.

The natural axiomatisation of this structure states that $j:\HH\to \C$ is a surjection such that for every $g_1,\ldots,g_n\in \mathrm{GL}^+(2,\Q)$ there is an irreducible  algebraic curve $C_{g_1\ldots g_n}\subset \C^{n+1}$ over  $\Q(S),$ where $S$ is the set of special points, such that  
$$\la y_0,y_1,\ldots,y_n\ra\in C_ {g_1\ldots g_n}\Leftrightarrow \exists  \tau\in \HH\ y_0= j(\tau), y_1= j(g_1\tau), \ldots y_n=j(g_n\tau).$$

This is given by a list of first-order sentences. Finally, the $L_{\omega_1,\omega}$-sentence 
$$j(\tau)=j(\tau') \Leftrightarrow \bigvee_{g\in \mathrm{SL}(2,\Z)} \tau'=g\tau$$
states that the fibres of $j$ are  $\mathrm{SL}(2,\Z)$-orbits.

Adam Harris proves in \cite{AH1} an analogue of Theorems A above.
\epk
\bpk \label{Aj} {\bf Theorem A$j$}  {\em The natural axiomatisation of 
$(\HH, j, \C_{\rm field}),$ the two-sorted structure for the $j$-invariant, defines an uncountably categorical AEC.}

The structure of the proof is similar to the proofs discussed above. The key model-theoretic tool is Theorem B$^*$. The appropriate thumbtack lemma takes the following form. 

{\bf Theorem} {\em Let $\A$ be an abelian variety defined over $\kk,$ a finitely generated extension
of $\Q$ or a finitely generated extension of an algebraically closed field $L$, such that $\A$ is
a product of $r$ non-isogenous elliptic curves (with $j$-invariants which are transcendental
over $\kk$ in the second case). Then the image of the Galois representation on the Tate
module $T(\A)$ is open in the Hodge group $\mathrm{Hg}(\A)(\hat\Z).$}

The {\em Hodge group} is a subgroup of the {\em Mumford-Tate group}, for definitions see  e.g. \cite{De}.

This theorem, essentially a version of the {\em adellic Mumford-Tate conjecture,} is explained in  \cite{AH1} as a direct consequence of a version of the hard theorem of Serre mentioned above and a further work by Ribet.  Consequently Harris deduces the categoricity theorem A$j.$

What is even more striking, that assuming the statement of  theorem A$j$ does hold, \cite{AH1} deduces the statement of the  adellic Mumford-Tate conjecture as a consequence. This sort of equivalence was observed in \cite{Zrav} for categoricity statements  for semi-Abelian varieties. The case considered by Harris, the simplest case of a Shimura variety, looks very different. And yet a similar tight connection between the model theory and the diophantine geometry of the $j$-invariant is valid. 
 
\epk
\bpk Most recent paper \cite{DH} by C.Daw and  A.Harris extends Theorem A$j$
to much broader class of two sorted structures, replacing  $ \C_{\rm field}$ by
an arbitrary {\em modular curve} and $j$ by the corresponding modular function.

Moreover,  C.Daw and  A.Harris extend the formalism of the two-sorted structure to any Shimura variety and show that the respective categoricity statement is equivalent to an (unknown) open image condition on  certain Galois representations. This equivalence
between model-theoretic statements and those from arithmetic may be seen
as a form of justification for the latter.
\epk
\bpk We would like to remark that in terms of the discussion in \ref{d1}, \ref{d2} the axiomatisation \ref{jaxiom}  assumed in Theorem A$j$ is not the most natural one. Its language includes constants naming special points, as a result of which the action of $\mathrm{Gal}(\tilde{\Q}:\Q)$ on special points is not seen in the automorphism group of the structure. Working in a more basic language one would require a stronger version of the corresponding ``Thumbtack Lemma'' and thus a deeper statement of Hodge theory.
 
\epk
\bpk A natural question (also asked by the anonymous referee) arises at this point: what would the analogue of  Theorem \ref{main1} (one-sorted case) be in the case of $j$ in place of $\exp,$ or even in in the context of \ref{d0}, for elliptic curves and Abelian varieties?

The case of an elliptic curve is expected to be quite similar to \ref{main1}. One needs just to try and prove the categoricity statement for the structure $\left(\C;,+,\cdot, \mathfrak{P}(\tau, z) \right),$ where $\tau$ is the modular parameter of the elliptic curve in question and
$ \mathfrak{P}(\tau, z) $ is the Weierstrass function in variable $z.$  The abelian
variety case shouldn't be very different; just replace $ \mathfrak{P}(\tau, z) $ by the $g$ functions $p_1(\bar z),\ldots, p_g(\bar z)$ of variables $\bar z$ on $\C^g,$ $g$ the genus of the variety, which make up the exponentiation map from $\C^g$ onto the variety.  The necessary form of the generalised Schanuel condition is known (but not the proof of it!), see the discussion in \ref{discuss}. 

In case of $j$ there is an essential difference: the natural domain of definition of 
$j$ is $\HH,$ the upper halfplane of the complex numbers. Since $(\C; +,\cdot, \HH)$ interprets the reals,  $(\C; +,\cdot, j)$ can not be excellent. Hence one needs to look for an extension of  $j$ to the complex plane minus a countable set in a way that preserves its modularity and the corresponding generalised Schanuel
condition (see \ref{discuss}).     This would be very desirable  for number theoretic reasons,  see Yu.Manin's paper \cite{Man} where a related problem is discussed.
An easy argument shows that this is impossible.  

Still, we may speculate, a possible solution to the problem would be to consider a ``multivalued function'' $J$ instead of $j.$ More precisely, define $J(x,y)$ to be a binary relation on $\C$ such that 

 for $x\in \HH,\ y\in \C$ $$J(x,y)\mbox{ iff }\exists g\in \mathrm{GL}^+_2(\Q)\  y=j(gx).$$
We want to extend $J$ to $\C$ so that

 for all $x\in \C\setminus \Q,\ y\in \C$ and   $g\in \mathrm{GL}_2(\Q),$
$$J(x,y)\to  J(gx,y),$$

and  for all $x\in \C\setminus \Q,\ y_0,y_1\in \C$ 
$$J(x,y_0)\ \& \ J(x,y_1)\to \exists g\in \mathrm{GL}^+_2(\Q)\ \la y_0,y_1\ra\in C_g,$$
where $C_g$ is as in \ref{jaxiom}.

For such a $J$  the same generalised Schanuel condition $\mathrm{SCH}_J$ as
for $j$ makes sense. Now it is possible to formulate  $\mathrm{EC}_J$ and have the analogue of \ref{main1} as the conjecture for the one-sorted structure   $(\C; +,\cdot, J).$  

The evidence supporting such a conjecture can be seen in the recently (formulated and) proved  {\em Ax - Schanuel} statement for the $j$ invariant in the differential field setting, \cite{Pila}, Conjecture 7.12, the proof by J.Tsimmermann and J.Pila to appear elsewhere.  

The tight relationship between the statement of \ref{main1}, the theorem of Ax (the Ax - Schanuel statement) and the model theory of differential fields has been given a deep analysis in the work of J.Kirby \cite{Ki1}. It should be also applicable in the case in question.  

\epk
\bpk Finally, we would like to discuss the first-order versus non-first-order (AEC) alternative in choosing a formalism to develop the model theory of pseudo-analytic  structures as above. 

In this regard there is a substantial difference between results of type A (Theorems A) and the main result about the ``one-sorted'' pseudo-exponentiation stated in \ref{main1}. In the first situation, as shown in \cite{BaHaPi}, one can use essentially the same techniques to classify models of the first order theory of the two-sorted structure in question.

In the one-sorted case the key property on which the whole model-theoretic study relies is  Hrushovski's
inequality (recall \ref{Hr}) or in more concrete form the relevant Schanuel's condition.  
The ``decomposition'' approach as in \ref{d2} can still  be attempted but, as was noted already in \cite{ZParis}, to formulate Schanuel's condition in the first order way one  requires certain degree of diophantine uniformity, which was formulated in \cite{ZLMS}   as the Conjecture on Intersections in Tori, CIT. \cite{ZParis} discusses a broader formulation of this which includes semi-Abelian varieties. An equivalent to CIT conjecture was later formulated by E.Bombieri, D.Masser and U.Zannier,  and a very general form, which covers the whole class of mixed Shimura varieties, formulated by R.Pink. Without giving precise definition of {\em special subvarieties} (the second part of this paper is devoted to this, see section \ref{ss}) we formulate what is currently referred to as the Zilber-Pink conjecture.
\epk
\bpk \label{ZP}{\bf Conjecture Z-P.} {\em Let $\X$ be an algebraic variety for which the notion of special subvarieties is well-defined. 
For any algebraic subvariety $V\subs \X^n$ there is a finite list of special subvarieties
$S_1,\ldots, S_m\subsetneq \X^n$
such that, given an arbitrary special subvariety $T\subset \X^n$ and an irreducible component $W$ of the intersection $V\cap T,$ either $\dim W=\dim V+\dim T-\dim \X^n$ (a typical case), or 
 $W\subs S_i$ for some $i=1,\ldots,m$ (in the atypical case $\dim W>\dim V+\dim T-\dim \X^n$).}
 
 This is a fundamental diophantine conjecture, in particular containing the Mordell-Lang and Andr\'e - Oort conjecture. 
 
 The analysis of CIT in \cite{ZLMS} specifically concentrates on $\X=\mathbb{G}_m,$ the algebraic torus (and characteristic 0), in which case the special subvarieties $S$ of   $\mathbb{G}_m^n$ are  subvarieties of the form
 $T\cdot t,$ where $T\subs \mathbb{G}_m^n$ is a subtorus and $t$ a torsion point in  $\mathbb{G}_m^n.$

\epk
\bpk In \cite{KZ} J.Kirby and the present author study the first order theory of the field with pseudo-exponentiation, $\F_{\exp},$ whose non-elementary theory is described in \ref{main1}. In this case, due to the richer language the structure on the kernel of $\exp$ is highly complex, more precisely, in the standard model it is effectively the ring of integers. So the theory of the kernel is the {\em complete arithmetic}.
Nevertheless we can aim at ``decomposing'' $\mathrm{Th}(\F_{\exp})$ into complete arithmetic and the theory
{\em modulo the kernel}. The main result states that assuming CIT such a description is possible. Moreover,
the theory modulo the kernel is quite tame and even superstable in a certain sense. The further analysis
shows that the key property of the structure, Schanuel's condition, is first order axiomatisable over the kernel  if and only if CIT holds.

The same must be true for other one-sorted pseudo-analytic structures, that is the study of the first order theories of such structures depend on the corresponding generalisations of CIT. On the other hand, the result of \cite{KZ} shows that any progress in the studies of the first order theories should shed light on the corresponding Z-P conjecture.
\epk

\bpk Finally, we conclude this section with the remark that to study the model theory of a  one-sorted pseudo-analytic structure one needs to know at the very least the statement of the corresponding Schanuel's condition. This is true both for AEC and the first order settings. A closer look at this problem quickly relates it with the problem of defining the notion of special subvarieties. The second part of the paper deals exactly with both issues. 
\epk

\bpk {\bf Introduction  to the second part}.

Consider a complex algebraic variety $\X$ and a semi-algebraic\footnote{Semi-algebraic sets are
by definition subsets of $\R^n$ that can be represented as Boolean combinations of subsets defined by the inequality $p(x_1,\ldots, x_n)\ge 0$ for $p$ a polynomial over $\R.$} set $\U$
which is also an open subset of $\C^m,$ some $m.$  $\X$ will be treated as a structure with the universe $\X$ and $n$-ary relations given by Zariski closed subsets of $\X^n.$  This is a classical Noetherian Zariski structure in the sense of \cite{Zbook}.
The structure on
$\U$ we define in a more delicate way so that eventually
 it is  Zariski  of analytic type.

Our full structure consists of two sorts $\U$ and $\X$ with an analytic surjective mapping 
$\Fp: \U\to \X$ connecting the sorts in a ``nice'' way so that $(\X,\U,\Fp)$ is model-theoretically tame
(in general $\Fp$ should be a correspondence, but in this text we deal with a mapping only).

By analogy with the theory of mixed Shimura varieties
call an irreducible Zariski closed $S\subset \X^n$ { weakly special}  if 
there is an analytically irreducible $\check{S}\subset \U^n$ such that $\Fp(\check{S})=S$ and $\check{S}$ is semi-algebraic (equivalently, definable in $\U$). Analogously but with more work one defines  {\em special} subsets $S\subseteq \X^n.$ The counterparts $\check{S}\subseteq \U^n$ of special sets $S$ are called {\em co-special}. Assuming ``niceness'' of the definitions involved the co-special subsets form a Zariski geometry (see Theorem~\ref{fcor2} below), so satisfy the Trichotomy principle:  its strongly minimal subsets $Y$ can be classified as being of one of the three types:
\begin{itemize}
\item trivial type: essentially a homogeneous space of a countable group;
\item linear type: typically a commutative group;
\item algebro-geometric type: $Y$ is a complex algebraic curve with  $n$-ary relations on $Y$ given by Zariski closed subsets of $Y^n.$  
\end{itemize}

 In particular, we have a well-defined notion of a combinatorial dimension (rank) $d_{\rm Spec}(u_1,\ldots,u_n)$ of a tuple $\la u_1,\ldots,u_n\ra\in \U^n.$ We also have
a combinatorial dimension (transcendence degree) $\trd_\Q(y_1,\ldots,y_n)$ for tuples $\la y_1,\ldots,y_n\ra\in \X^n$ as well as for  tuples $\la y_1,\ldots,y_n\ra\in \U^n.$

From model-theoretic point of view $\Fp$ can be interpreted as a ''new`` relation on an algebraically closed field $(\C,+,\cdot)$ and in order 
for the structure $\C_\Fp:=(\C,+,\cdot, \Fp)$ to be model-theoretically ``nice''  the only construction known today (after more than 20 years since
\cite{Hr0}) is to employ the Hrushovski fusion method. In  our context  Hrushovski's 
construction suggests the following.

Introduce a relevant Hrushovski predimension. For $u_1,\ldots,u_n\in \U$
set 
$$\delta(u_1,\ldots,u_n)=\trd_\Q(u_1,\ldots,u_n, \Fp(u_1),\ldots,\Fp(u_n))-d_{\rm Spec}(u_1,\ldots,u_n).$$

The standard assumption (Hrushovski's inequality) then would be \be \label{F-Sch} \delta(\bar u)\ge 0,\mbox{ for all }\bar u\subset \U\ee
which is the {\em generalised Schanuel conjecture} (in the sense of \cite{ZParis}) corresponding to our $\Fp.$
\epk
\bpk\label{discuss}  We can rewrite (\ref{F-Sch}) as
\be \label{F-And}\trd_\Q(u_1,\ldots,u_n, \Fp(u_1),\ldots,\Fp(u_n))\ge d_{\rm Spec}(u_1,\ldots,u_n)\ee
and in this form compare it with Andr\'e's generalisation of Grothendieck period conjecture (\cite{An}, 23.4.1)
\be \label{Andre} \trd_{\Q}k(\mathrm{periods}(\X^n))\ge \dim G_{\mathrm{mot}}(\X^n)\ee
where $\trd_{\Q}k(\mathrm{periods}(\X^n))$ is the transcendence degree of periods of $\X^n$ with parameters in $k$ and
$G_{\mathrm{mot}}(\X^n)$ is the motivic Galois group of $\X^n.$

C. Bertolin in \cite{Be} considers the  1-{\em motives} case of  Andr\'e's conjecture and translates it in the form  that  generalises Schanuel's conjecture  covering the case of elliptic functions in combination with exponentiation. Andr\'e discusses this in \cite{An}, 23.4.2, and later on in 23.4.4 shows that 
his conjecture in case of Shimura variety $\mathit{A}_g$ (the moduli space of Abelian varieties of dimension $g$) and  the corresponding invariant $j_g$ defined on  the Siegel half-space $\HH_g,$ one has
$$\trd_\Q(\tau, j_g(\tau))\ge \dim_\C(j_g(D_\tau)),$$
where $D_\tau$ is the smallest co-special subvariety of $\HH_g$ which contains $\tau.$ These special cases of Andr\'e's conjecture completely agree with our (\ref{F-And}).    

This comparison suggests that motivic objects can be explained in terms of the geometry of co-special sets and co-special points. An approach to the classification of these geometries is the subject of this paper, see Theorems \ref{fcor2} and \ref{fcor3}. In particular we show that any such geometry is a combination of the three basic types of geometries of the Trichotomy principle.  
\epk

\section{Special coverings of algebraic varieties}

\bpk \label{f00} {\bf The setting.}

 We consider $\C$  both in complex and in real co-ordinates, so we can define semi-algebraic subsets and relations in $\C.$ 
\begin{itemize}
\item[A.]
\begin{itemize}
\item[(i)]
 We are given an open 
 $\U\subs  \C^m$ and with the complex structure induced from $\C$  the set
 $\U$ can be considered a complex manifold,  and the same $\U$  we view in real co-ordinates.
\item[(ii)]
We are also given  a smooth complex algebraic variety $\X$ and an analytic 
surjection   $\Fp:\U\to \X$ with discrete fibres.
\item[(iii)]
We assume that a  group $\Gamma$ acts on $\U$ by biholomorphic transformations and discontinuously. The fibres of $\Fp$ are orbits of $\Gamma.$    
  \end{itemize}
\item[B.] 
Further on we assume the existence of a  {\em semi-algebraic fundamental domain} $\V\subs \U$  such that 
\begin{itemize}
\item[(i)]
$\dot{\V}\subs \V\subs \overline{\V},$ where $\dot{\V}$ is the interior of $\V$ and  $\overline{\V}$ is the closure of  $\dot{\V}$ in the metric topology. We also assume that $\overline{\V}$ is semi-algebraic.

\item[(ii)] For each $\gamma\in \Gamma$ the set ${\V}_\gamma=\gamma\cdot \V$ and the
 restriction of the map $z\mapsto \gamma\cdot z$ on 
$\overline{\V}$ are semi-algebraic,
$$\dot{\V}_\gamma\cap \dot{\V}=\emptyset,\mbox{ for }\gamma\neq 1$$ 
 and
$$\U=\bigcup_{\gamma\in \Gamma} {\F}_\gamma.$$


 \end{itemize}

\item[C.]  The set
$$\Delta:=\{ \gamma\in \Gamma: \overline{\V}\cap \overline{\V}_\gamma\neq \emptyset\}\mbox{ is finite}.$$
\end{itemize}
\epk
\bpk \label{Example-shimura} {\bf Examples. } A large class of examples is the class of {\em arithmetic varieties} $\X:=\Gamma\backslash \U,$ where $\U$ is a  Hermitian symmetric domain and $\Gamma$ is an arithmetic subgroup of the real adjoint group Lie  of biholomorphisms of $\U.$ It would be impractical to explain these classical notions in the current paper. We just refer to \cite{Borel} for the introduction to the subject. 

Arithmetic varieties have fundamental domains in the form of
  {\em Siegel sets}, which are semi-algebraic. Moreover, the condition \ref{f00}.C is satisfied. See \cite{Borel}, Thm 13.1.

This class includes all {\em Shimura varieties} $\X.$  
\epk


  Do mixed Shimura varieties satisfy A-C? 

\bpk \label{Example1}
{\bf Example. }
 Let $$\exp:\mathfrak{sl}(2,\C)\to \text{SL}(2,\C)$$
be the Lie exponentiation. The restriction of $\exp$ to the nilpotent part of $\mathfrak{sl}(2,\C)$ is
algebraic, in fact it is a map $z\mapsto 1+z$ on a nilpotent matrix $z$ ($1$ is the unit matrix).

Let $N\subset   \text{SL}(2,\C)$ be the set of unipotent elements of  $\text{SL}(2,\C)$ (including $1$).
We define $$\X=\text{SL}(2,\C)\setminus N\mbox{ and }\U=\mathfrak{sl}(2,\C)\setminus\mathrm{Ln}(N)$$
Clearly $\exp(\U)=\X$ and $\exp$ is unramified on $\U.$

We want to define $\Gamma$ and a  fundamental domain for  $\exp: \U\to \X.$ Looking at the  Jordan normal form we see that an arbitrary $a\in \mathfrak{sl}(2,\C)$
is  of the form $$a=\left( \begin{array}{ll} x\ \ \ 0\\ 0 \ -x\end{array}\right)^g, \ \ g\in \text{SL}(2,\C)$$

We write  $b^{\text{SL}(2,\C)}$ for the conjugacy class $\{ b^g:g\in  \text{SL}(2,\C)\}.$

Define 
$$\V=\left\{ \left( \begin{array}{ll} x\ \ \ 0\\ 0 \ -x\end{array}\right)^{ \text{SL}(2,\C)}: x\in \C, \  x\neq 0\mbox{ and } 0\le \mathrm{Im} x<2\pi\right\}$$

The uniqueness, up to the order of eigenvalues, of diagonalisation implies
that  $\exp$ is injective on $\V.$ 

Define $\gamma_k:\U\to \U$ by setting
$$\gamma_{k}: \left( \begin{array}{ll} x\ \ \ 0\\ 0 \ -x\end{array}\right)^g\mapsto
\left( \begin{array}{ll} x+2\pi ik\ \ \ \ 0\\ 0 \ \ -x-2\pi ik\end{array}\right)^g.$$

The maps are well-defined since the cetralisers of  $\left( \begin{array}{ll} x\ \ \ 0\\ 0 \ -x\end{array}\right)$ and $
\left( \begin{array}{ll} x+2\pi ik\ \ \ \ 0\\ 0 \ \ -x-2\pi ik\end{array}\right)$ coincide.
 The restriction of $\gamma_k$ on $\V$ is a semi-algebraic
 since determining eigenvalues of an element of a matrix group is a definable operation over the field of reals. 

 The set $\Delta$  of \ref{f00}C is equal to $\{ \gamma_1, 1, \gamma_{-1}\}.$ 

\medskip

{\bf Remark.} In this example $\U$ is not semi-algebraic since the subset of diagonal matrices in $\U$ (which can be defined by an algebraic formula) is in algebraic bijection with $\C\setminus 2\pi i\Z.$

\epk
\bpk \label{dim} We are going to use the notion of dimension $\dim_\R T$ for certain ``nice''  subsets of $T$ of $\U,$ $\X$ and their Cartesian powers. Obviously, this notion is applicable when $T$ is semi-algebraic. More generally, note that given an open ball  $B\subset \U^n$  of small enough radius, the two sorted structure
$(B, \Fp(B))$ along with the map $\Fp$ restricted to $B$ is definable in the o-minimal structure $\R_{an}$
(here and below we use the same notation for
the map $\la u_1,\ldots,u_n\ra\mapsto  \la \Fp(u_1),\ldots,\Fp(u_n)\ra$ for all $n\ge 1)$).     

Now we can apply the dimension $\dim_\R$ in the sense of  $\R_{an}$ to any subset definable in this structure. In particular, if $T$ is a semi-algebraic subset of $\U^n$ (in respect to the embedding 
$\U\subs \C^m$) then  $\dim_\R(T\cap B)$ and $\dim_\R(\Fp(T\cap B))$ are well defined. We can now define
$$\dim_\R \Fp(T)=\max_B \dim_\R\Fp (T\cap B).$$ 

\medskip

A direct consequence of discreteness of fibres is the following.

\medskip

{\bf Fact.} {\em The map $\Fp$ preserves dimension. More precisely, 
for every semi-algebraic $T\subs \U^n,$
$\dim_\R \Fp(T)=\dim_\R T.$

For every semi-algebraic $S\subs \X^n,$
$\dim_\R \Fp\inv(S)=\dim_\R S.$}
\epk

\bpk \label{glu} {\bf Lemma.} {\em  {\rm (i)} The equivalence relation 
$\sim$ on $\overline{\V}$ defined as 
$$ u_1\sim u_2\mbox{ iff } \Fp(u_1)=\Fp(u_2)$$
 is semi-algebraic.

{\rm (ii)} There is a (non-Hausdorff) topology $\mathfrak{T}$ on $\overline{\V}$
\begin{itemize}
\item
the base open sets of $\mathfrak{T}$ are of the form $\overline{\V}\cap \Gamma\cdot B$ for $B\subs \U$ an open ball;
\item  $\mathfrak{T}$ induces a Hausdorff topology $\mathfrak{T}/_\sim$ 
on the quotient $\overline{\V}/_\sim$ with the base of open sets of the form
$(\overline{\V}\cap \Gamma\cdot B)/_\sim$ and
 the canonical map associated with $\sim,$  
$$\tilde\Fp:\overline{\V}/_\sim\to \X,$$  is a homeomorphism with respect to the metric topology on $\X;$
\item the restriction of $\Fp$ to $\overline{\V},$
$$\Fp: \overline{\V}\to \X,$$
is an open and closed map.
\end{itemize}

 {\rm (iii)} There are a semi-algebraic set $\tilde{\V},$ $\dot{\V}\subs\tilde{\V}\subs \overline{\V},$ and a bijective semi-algebraic correspondence  
$\mathbf{i}:\overline{\V}/_\sim\to \tilde{\V} $ such that $\tilde{\Fp}=\Fp\circ \mathbf{i}.$  

}
 
 {\bf Proof.} Recall $\Delta=\{ \gamma\in \Gamma: \gamma\overline \V\cap \overline{\V}\neq \emptyset\}.$
 
 We can redefine equivalently for $u_1,u_2\in \overline{\V}$
 $$u_1\sim u_2\mbox{ iff } \exists \gamma\in \Delta\ u_2=\gamma\cdot u_1.$$
Since  $\Delta$ is finite the relation is semi-algebraic. This proves (i).

\medskip

(ii) It is clear that $\mathfrak{T}$ is a weakening of the metric topology on $ \overline{\V}$ since $\Gamma B$ is open in $\U.$

Let $\Fp(u)=x,$ $\Fp(v)=y$ and $V_x,$ $V_y$ open non-intersecting neighbourhoods in $\X$ around $x$ and $y$ correspondingly. 

Then  we can find balls
$B_u\subs \Fp\inv(V_x)$ and $B_v\subs \Fp\inv(V_y)$ around $u$ and $v$ correspondingly. Since
$\Fp\inv(V_x)\cap \Fp\inv(V_y)=\emptyset$ and both are invariant under the action of $\Gamma,$ we have $\Gamma B_u\cap \Gamma B_v=\emptyset.$ 
This proves that $\mathfrak{T}/_\sim$ is Hausdorff.

The argument above also shows that the inverse image of an open subset of $\X$ under $\Fp$ is $\mathfrak{T}$-open in $\overline{\V},$ and  the inverse image of an open subset of $\X$ under $\tilde{\Fp}$ is $\mathfrak{T}/_\sim$-open in $\overline{\V}/_\sim.$ This prove that $\Fp$ and $\tilde{\Fp}$ are continuous maps in the relevant topologies.

We can also characterise $\mathfrak{T}$ in terms of its base of closed subsets. We can take for these the sets of the form  $\overline{\V}\cap (\U\setminus \Gamma B)$ for $B$ a ball in $\U.$  By the above $\Fp$ sends such a closed set onto a closed set, and the same true for $\tilde{\Fp}.$

It follows that $\tilde{\Fp}: \overline{\V}/_\sim\to \X$ is a homeomorphism.

\medskip

(iii) Since for any semi-algebraic surjection $f:X\to Y$ there is a semi-algebraic section $s: Y\to X$  
 there exists a semi-algebraic section $\mathbf{i}:\overline{\V}/_\sim\to  \overline{\V}$ inverse to the quotient-map.
 Set  $\tilde{\V}:=\mathbf{i}(\overline{\V}/_\sim).$ $\Box$
\epk
\bpk \label{remT} {\bf Remark.} In particular, we have seen in the proof of (ii) that given a closed subset $C\subs \U$ invariant under the action of $\Gamma$ (that is $\Gamma C=C$) the image $\Fp(\overline{\V}\cap C)$ is closed in $\X.$
 
\epk

\bpk \label{newU} Note that $\gamma\cdot\tilde{\V}\cap\tilde{\V}=\emptyset$ for $\gamma\neq 1,$ since
$\sim$ is trivial on $\tilde{\V}$ by definition. Without loss of generality {\bf we assume from now on} that
$ \tilde{\V}=\V$ and thus have a stronger condition for \ref{f00}.B(ii): 
$$\gamma\V \cap \V=\emptyset\mbox{ for }\gamma\neq 1.$$
\epk

\section{The geometry of weakly special subsets}

\bpk {\bf Definition.} We say $S\subs \X^n$ is weakly special, if $S$ is Zariski closed irreducible and $\Fp\inv(S)\cap {\V}^n$ is semi-algebraic.

Note that if $\Fp\inv(S)\cap \gamma{\V}^n$ is semi-algebraic for $\gamma=1$, it is semi-algebraic for every $\gamma\in \Gamma^n.$ It is immediate by \ref{f00}.B(ii).
\epk

\bpk \label{diag}
Note that every point and the diagonal of $\X^2$ are weakly special.

\epk
\bpk \label{cells}
We use the definition of a {\em cell} $C\subs \R^m$ and the theorem of semi-algebraic cell decomposition
in the context of the field of reals, $(\R,+,\cdot,\le),$ see \cite{vdD}. Such a $C$ can be represented as an intersection $X\cap U,$ where $U\subs \R^m$ is semi-algebraic open and $X$ the set of real points of a Zariski closed set. We may also assume that $C$ is irreducible as  real analytic variety (that is can not be presented as a proper union of two analytic subvarieties).

\epk
\bpk \label{f01} {\bf Lemma.} {\em Let $C\subs  \U^n$  be a 
cell, $X$ a real analytic subset of $\U^n$ and $\dim_R(X\cap C)=\dim_R(C).$   Then $C\subs X.$}

{\bf Proof.}
Immediate from the fact that $C$ is irreducible real analytic. $\Box$

\epk
\bpk \label{XG} {\bf Lemma.} {\em Let $X$ be a real analytic subset of an open $G,$ both semialgebraic. Let $X=\bigcup_iX_i$ be its decomposition into irreducible components. Then there are only finitely many non-empty $X_i$ and all of them are semi-algebraic.}

{\bf Proof.} Note that the set $X^{\rm sing}$ of singular points is definable, i.e. semi-algebraic. 
Let $X'= X\setminus X^{\rm sing}$ and $G'=G\setminus X^{\rm sing},$ which is also open. Consider a cell decomposition $$X'=\bigcup_{j=1}^m C_j$$ 
and let $X'_i=X_i\cap G'.$

For each $C_j$
there must be an $X_{i_j}$ such that $\dim_R(X_i\cap C_j)=\dim C_j.$ By \ref{f01} $C_j \subs  X_{i_j}.$ This proves $X'= \bigcup_{j=1}^mX_{i_j}.$

Since non-empty intersections $X_i\cap X_k$ for $i\neq k$ are subsets of  $X^{\rm sing},$   the  decomposition of $X'$ into the $X'_i$ is  disjoint.
Then $X'_{i_j}$ either does not intersect $C_l$ or contains it and is equal to $X'_{i_l}.$ It follows, that
$X'_{i_j}$ is a union of finitely many cells, so is semi-algebraic. 

It remains to note that the metric closure of  $X'_i=X_i\setminus X^{\rm sing}$ is equal to $X_i$ and the $X'_i$ are semi-algebraic, the $X_i$ are semi-algebraic. $\Box$
\epk

\bpk \label{R} {\bf Lemma.} {\em Suppose $R\subs \X^n$ is Zariski closed and $\Fp\inv(R)\cap \V^n$ semi-algebraic. Let $\Fp\inv(R)=\bigcup_mT_m$ be the decomposition of the analytic subset $\Fp\inv(R)$ of $\U^n$ into irreducible components. Then $T_m\cap \V^n$ is semi-algebraic for every $m$ and there are only finitely many of such subsets. 

Moreover, given  $$R=R_1\cup\ldots\cup R_k,$$
the decomposition into Zariski irreducible components, the sets $\Fp\inv(R_i)\cap \V^n$ are 
semi-algebraic.
}

{\bf Proof.} Let $$\Fp\inv(R)\cap \overline{\V}^n=C_1\cup \ldots\cup C_p$$ be a decomposition into semi-algebraic cells. 



Suppose $T_{m}\cap \overline{\V}^n\neq \emptyset.$ Then there is a cell, say $C_1$ of the decomposition such that $\dim_\R(T_{m}\cap C_1)=\dim_R(T_{m}\cap \overline{\V}^n).$ 

Let 
$G$ be an open semi-algebraic subset of $\U$ such that $C_1$ is a real analytic irreducible subset of $G.$ Then by \ref{f01} $C_1\subs  T_m\cap G.$ Since $T_m\cap  \overline{\V}^n$  is closed
the closure $\overline{C}_1\subs  T_m\cap G.$ But since  $C_1\cup \ldots\cup C_p$ is closed in $\U$,
we have that  $\overline{C}_1$ is a union of cells in this decomposition (boundary cells), say $\overline{C}_1=\bigcup_i\le jC_i$ for some $j\le p.$ Now, if $T_{m}\cap \overline{\V}^n\neq \emptyset$  we  consider $T'_m:=T_{m}\setminus \overline{C}_1$ an analytic subset of the open set $\U'= \U\setminus \overline{C}_1.$ Note that the cells are disjoint, so $C_{j+1},\ldots,C_p\subs {\U'}^n\cap \overline{\V}^n.$

We can again  find a cell, say $C_{j+1}$ such that
 $\dim_\R(T_{m}\cap C_{j+1})=\dim_R(T'_{m}\cap \overline{\V}^n)$ and repeating the argument will get
 $C_{j+1}\subs T'_m.$ After finite number of steps we get that
$T_{m}\cap \overline{\V}^n= \bigcup_i\le kC_i,$ for some $k\le p,$ up to renumeration of the decomposition. Claim proved.

It follows that every $T_{m}\cap \overline{\V}^n$ is semi-algebraic and there are only finitely many such.

Now note that by uniqueness of irreducible decomposition the analytic set
$\Fp\inv(R_i)$ is a union of some subfamily of the $T_m$'s. 
Thus  $\Fp\inv(R_i)\cap {\V}^n$ are  semi-algebraic. $\Box$
\epk
\bpk \label{R2} {\bf Lemma.} {\em Let $R\subs \X^n$ be a constructible subset (a Boolean combination of Zariski closed sets) and $T:=\Fp\inv(R)\cap \V^n$ be semi-algebraic. Then $R$ is a Boolean combination of weakly special subsets of $\X^n.$}

{\bf Proof.} $R$ can be represented as a finite union of constructible sets of the form
$R_i\setminus P_i,$ where $R_i$  Zariski closed irreducible and $P_i$ are Zariski closed subset of $R_i,$
$\dim_\C(R_i)<\dim_\C(R_i).$

Consider the metric closure $\overline{R}.$ By general facts  
$\overline{R}=\bigcup_{i\le k}R_i.$ 

On the other hand the $\mathfrak{T}$-closure $\overline{T}$ of $T$ is semi-algebraic and the two are homeomorphic by ${\Fp}:  \overline{T}\to \overline{R}.$
It follows by \ref{R} that the components $R_i$ are weakly special and $T_i:=\Fp\inv(R_i)\cap \V^n$ are semi-algebraic. 

By homeomorphism,  $ \bigcup_{i\le k}T_i$ is the closure of $T.$ Moreover, 
the semi-algebraic set $T_i\cap T$ corresponds to  $R_i\setminus P_i,$ that is
$T_i\cap T=T_i\setminus Q_i$ for some  $Q_i\subset T_i,$ $\Fp(Q_i)=P_i.$ By construction $Q_i$ is semi-algebraic. Hence, by \ref{R}, $P_i$ is a union of weakly special sets. $\Box$
\epk

\bpk \label{f1} {\bf Proposition.} {\em 
Let $S,S_1$ and $S_2$ be weakly special subsets of $\X^n.$ Then

(i) the irreducible components of $S_1\cap S_2$ are weakly special.

(ii) given the projection $\pr: \X^n\to \X^{n-1}$ one can represent $$\pr S=R\setminus P,$$
where $R$ is weakly special  and $P$ is
 a Boolean combination  of weakly special sets, and  $\dim_\C P<\dim_\C R.$

(iii) Let $\la x_1,\ldots, x_{n-1}\ra$ be a generic point of $\pr S,$ $S_{x_1\ldots x_{n-1}}$  the fibre over $\la x_1\ldots x_{n-1}\ra$ and
$d=\dim S_{x_1\ldots x_{n-1}},$  the dimension of the generic fibre. Then
the set $$S^{(d)}=\{ \la y_1,\ldots,y_{n-1}\ra\in \pr S:\ \dim_\C S_{y_1\ldots y_{n-1}}= d\}$$
contains a subset  of the form $\pr S\setminus \bigcup_{i\le k}Q,$ where $Q_1,\ldots,Q_k$ are weakly special subsets of $\X^{n-1}.$
}

{\bf Proof.} (i) Since $\Fp\inv(S_1\cap S_2)\cap {\V}^n=\Fp\inv(S_1)\cap  \Fp\inv(S_2)\cap {\V}^n,$ the set $\Fp\inv(S_1\cap S_2) \cap {\V}^n$ is 
semi-algebraic. Given $S_1\cap S_2=P_1\cup\ldots\cup P_k,$
the decomposition into Zariski irreducible components we have by
 \ref{R} that $\Fp\inv(P_i)$ is semi-algebraic, for every $i\le k.$

\medskip

(ii) $\pr S$  is Zariski constructible by Tarski's theorem. On the other hand $\pr T$ is semi-algebraic, for $T=\Fp\inv(S)\cap \V^n.$ The bijection $\Fp$ on $\V$ commutes with $\pr,$ so $\pr T=\Fp\inv(\pr S)\cap \V^{n-1}.$

Now \ref{R2} proves that $\pr S$ is a Boolean combination of weakly special sets,  so a union of sets of the form $R_i\setminus P_i$ where $R_i$ is weakly special and $P_i$ is a union of weakly special. Since
$\pr S$ is irreducible, we get the desired.

\medskip

(iii)  The standard Fibre Dimension Theorem implies that  the set $S^{(d)}$ is
constructible and contains a Zariski open subset of the irreducible set $\pr S.$ On the other hand for  $T=\Fp\inv(S)\cap \V^n$ we can define the set
$$T^{(d)}=\{ \la y_1,\ldots,y_{n-1}\ra\in \pr T:\ \dim_\R T_{y_1\ldots y_{n-1}}= 2d\},$$
 which is in bijective correspondence to $S^{(d)}$ via $\Fp$ since $\Fp$ preserves the real dimension.
Now note that $T^{(d)}$ is semi-algebraic since dimension of fibres is definable in the field of reals too. 
Thus we get the conditions of \ref{R2} for  $S^{(d)},$ which proves that  $S^{(d)}$ can be presented as a Boolean combination of weakly special sets. The statement follows.$\Box$

\epk
\bpk {\bf Definition.} We denote ${\cal S}_w$ the family of weakly special sets.  $\X_{S_w}$ will stand for the structure
$(\X, {\cal S}_w)$ with the universe $\X$ and $n$-relations given by weakly special subsets of $\X^n.$

We consider  $\X_{S_w}$ a topological structure with dimension. The closed sets of the topology on $\X^n,$ any $n,$ are finite unions of weakly special subsets and the dimension that  of algebraic geometry.  

\medskip

$\V_{S_w}$ will stand for the structure
$(\V, {\cal S}_w)$ with the universe $\V$ and $n$-relations given by subsets $\Fp\inv(S)\cap \V^n$ for $S\in {\cal S}_w,$ weakly special subsets.

\epk

\bpk \label{frem}{\bf Remark.} 
 $\Fp$ restricted to $\V$ is an isomorphism between  $\V_{S_w}$ and $\X_{S_w}.$


\epk

\bpk \label{fcor2} {\bf Theorem.} {\em 

(i) $\X_{S_w}$ is a  Noetherian Zariski structure.  

(ii)  $\X_{S_w}$ has quantifier elimination and is saturated. 

}

{\bf Proof.} (i) We refer to \cite{Zbook} for the definition of a Zariski structure:   

The properties (L),  (DCC), (SI) and (AF) follow from Proposition~\ref{f1}(i) and \ref{diag} along with the obvious properties of Zariski topology in algebraic geometry as discussed in \cite{Zbook}, 3.4.1.

The property (SP) of projection is checked by \ref{f1}(ii). 
 
 
 The fibre condition (FC) is \ref{f1}(iii).

(ii) 
Every Zariski structure has QE. Saturation follows from the fact that  $\X_{S_w}$ is definable in the saturated structure, the field $\C.$

 $\Box$
  
\epk 
\bpk We call a weakly special $\Y\subseteq \X^n$  {\bf simple} if it has no proper infinite weakly special subsets. 

In model theoretic terminology this is equivalent, by \ref{fcor2}, to say that $\X_{S_w}$ is strongly minimal.
\epk

\bpk \label{simple}{\bf Lemma.} {\em Suppose 
$\X$ is {\bf simple}. 

Set  
$\X^{reg}_{S_w}\subseteq \X_{S_w}$ to be the substructure  on the subset
$\X\setminus \X^{\rm Sing}$ 
(by restrictions of predicates and constants of $\mathcal{S}_w$).

Then 
$\X^{reg}_{S_w}$ is a presmooth strongly minimal Noetherian Zariski structure. }


{\bf Proof.} For convenience of notations we rename below $\X^{reg}=\Y.$ 

Claim. Given 
 a weakly special subset $S$ of $\X^n$  either $\dim_\C S\cap \Y^n=\dim_\C S$ or $S\cap \Y^n=\emptyset.$
Moreover, if 
no projection $\mathrm{pr}_1: \X^n\to \X$
along $n-1$ co-ordinates   $\mathrm{pr}_1 S$ is a point of $\X^{\rm Sing},$ we get 
$\dim_\C S\cap \Y^n=\dim_\C S.$ 

Suppose towards a contradiction that $\dim_\C S\cap \Y^n<\dim_\C S.$ Then, since $S$ is 
irreducible,  $S$ is a subset of Zariski closed set $\X^n\setminus \Y^n,$ that is
$S\cap \Y^n=\emptyset,$ This is the same as to say that, up to the order of co-ordinates, 
$S\subseteq \Y^{n-1}\times \X^{\rm Sing}.$ Then the corresponding projection $\mathrm{pr}_1 S\subseteq 
\X^{\rm Sing}.$ Thus the Zariski closure $\overline{\mathrm{pr}_1 S}$ is a subset of the Zariski closed set
$\X^{\rm Sing}.$
But $\overline{\mathrm{pr}_1 S}$ is weakly special by \ref{f1}(ii), so by assumptions the embedding can only happen if $\mathrm{pr}_1 S$ is finite, so a point, which contradicts the assumption on $S.$ Claim proved.

We will call weakly special in regards to $\Y$ the sets  of the form $S\cap \Y^n$ for $S$ weakly special subsets of  $\X^n.$
Now we can check that the statements (i), (ii) and (iii) of \ref{f1} hold for $\Y.$  (i) and (ii) in regards to $\Y$ follow from (i) and (ii) for $\X$ by definition. (iii) follows by the Claim by the same argument as in the proof of \ref{f1}.

Finally, $\X^{reg}_{S_w},$ the structure on $\Y,$ is Noetherian Zariski by the same argument as in the proof of \ref{fcor2}. The structure is presmooth since the underlying quasiprojective algebraic variety $\X\setminus \X^\mathrm{Sing}$ is smooth, see \cite{Zbook}, 3.4.1.
 $\Box$
\epk

\bpk \label{fcor2iii} {\bf Theorem.} {\em The simple weakly special subsets  $Y\subs \X^n,$ all $n\ge 1,$ are classifiable as follows:  

 the geometry of the structure $\Y$ induced  from $\X_{S_w}$ is either trivial, or linear (locally modular),  or of algebro-geometric type,  
in which case
\begin{itemize}
\item[(i)]
 $\dim_\C Y=1,$ that is $Y$ is an algebraic curve, and 
\item[(ii)] 
 every irreducible Zariski closed subset of $Y^k,$ $k\in \N,$ is weakly special. 
\end{itemize}

}

{\bf Proof.} By the Weak Trichotomy Theorem (see \cite{Zbook}, Appendix, Thm B.1.43) the geometry of
the strongly minimal structure $\Y$ is either trivial, or locally modular, or a pseudo-plane is definable in $\Y.$ We need to analyse the latter case.

In this case by B.1.39 of \cite{Zbook} on a subset $S$ of $Y\times Y$ there is a definable family $L$ of ``curves'' $C_l,$ $l\in L,$ with Morley ranks $\Mrk Y=1,$ $\Mrk L=2$ and $\Mrk C_l=1$ for each $l\in L.$
By elimination of quantifiers in $\Y$ each $C_l$ is a Boolean combination of weakly special subsets, so we may assume
$C_l$ is of the form $C_l=R_l\setminus Q_l,$ where $R_l$ is weakly special and $Q_l$ a finite union of 
weakly special subsets. 

Note that once it is established that $\Mrk Y=1,$ we have an easy translation between the Morley rank and the (complex) dimension of definable subsets $Z$ of $\Y^n:$ \ $\dim_\C Z=\Mrk Z\cdot \dim_\C Y.$

 By a standard argument for generic $l\in L,$ we have  $\Mrk(\pr C_l)=1,$ for projections of $Y\times Y\to Y$ along the both co-ordinates. This condition is definable in the structure $\Y,$ hence we may assume 
$\Mrk(\pr C_l)=1$ for all $l\in L.$ In particular, the Zariski closure of $\pr C_l$ is $\Y$ for all $l\in L.$

 
Consider  the substructure $\Y^{reg}$ of $\Y$ obtained by removing the singular points of $Y.$         
By \ref{simple}  $\Y^{reg}$ is a presmooth strongly minimal Noetherian Zariski structure such that
the trace of the family $C_l: l\in L$ on the substructure gives  us a Morley rank 2 family of ``curves'' $\Mrk(C_l\cap \Y^{reg}\times \Y^{reg})=1$ This proves that $\Y^{reg}$ is not locally modular.

 The classification of  a non locally modular presmooth strongly minimal Zariski structure is given by
  \cite{Zbook}, Theorem 4.4.1 (originally \cite{HZ}). 
By the theorem
a structure $\F=(F,+,\cdot)$ of an algebraically closed field $\F$  is definable in $\Y^{reg}.$ But 
$\Y^{reg}$ by construction is definable in the field $\C$ of complex numbers. So $\F$ is definable in $(\C,+,\cdot).$ By standard model-theoretic fact we can definably identify $\F$ and $\C.$

The above implies that there is a rational function $r: \Y^{reg}\to \C$ 
such that any Zariski closed subset $R\subset \C^n$ 
is also definable in $\Y^{reg}$ 
By strong minimality of $\Y^{reg}$ the map $r$ must be finite (that is with finite fibres $r\inv(x)$ ).
This implies that $\dim_\C \Y^{reg}=\dim_\C \C=1,$ that is $\Y^{reg}$ is an a complex algebraic curve.
Recall that $\Y^{reg}=\Y\setminus \Y^{\rm Sing}$ and $\Y$ is an algebraic variety. It follows that $\Y$ is a complex algebraic curve. We have proved (i).

Moreover, it follows that $\Y^{\rm Sing}$ is a finite subset, so in terms of definability $\Y$ and $\Y^{reg}$ are equivalent. So the statement of 4.4.1 is applicable to $\Y.$

The Classification Theorem 4.4.1 also states that the field is ``purely definable", that is any subset of $\C^m$ which is definable in $\Y$ is definable
in the field $\C$ alone.

Pick an arbitrary Zariski closed subset $P\subs \Y^m.$ 
Then $r(P)\subs \C^m$ is definable
in $\C,$ so $r\inv (r(P))$ is definable in $\Y.$ By elimination of quantifiers $r\inv (r(P))$ is a Boolean combination
of weakly special subsets, hence
the Zariski closure
$\bar P$ of $r\inv (r(P))\subs \Y^m$ is a finite union of weakly special subsets, that is definable in $\Y.$ 

Obviously $P\subseteq \bar P,$ and indeed $P$ is an irreducible component of $\bar P.$ By  Lemma~\ref{f1}(i)  $P$ is 
weakly special. $\Box$
\epk

\section{Co-special geometry on $\U$}.
\bpk \label{setting} Given a weakly special $S\subset \X^n$ consider the analytic subset $\Fp\inv(S)\subs \U^n$ and its
decomposition into analytic irreducible components
$$\Fp\inv(S)=\bigcup_{i\in I_S}T_i.$$
Note that for all components $\dim_\C T_i\le \dim_\C S$ since $\dim_\C \Fp\inv(S)=\dim_\C S.$ 

We will call a component $T_i$ {\bf essential} if $\dim_\C T_i=\dim_\C S.$

We will call a component $T_i$ {\bf essential in $\gamma\V^n$ } if $\dim_\C T_i\cap \gamma\V^n=\dim_\C S.$

Note that since $T_i\subs \bigcup_{\gamma\in \Gamma^n}\gamma\V^n,$ every essential $T_i$ is essential in some
$\gamma\V^n,$ $\gamma\in \Gamma^n.$

By \ref{R} for a given $\gamma$ there is only finitely many $T_i$ essential in $\gamma\V^n.$

\epk
\bpk \label{B} {\bf Lemma.} {\em  In notation of \ref{setting} let $T_1,\ldots,T_k$ be the essential components of $\Fp\inv(S)$ intersecting $\V^n.$ Then for any $i\le k$ there is $\gamma_i\in \Gamma^n$ such
that $T_i=\gamma_iT_1.$

Every components of $\Fp\inv(S)$ intersecting $\V^n$ is essential. 
}

{\bf Proof.} Note that for each $i$ the set  $\Gamma\cdot T_i$ is analytic and so closed in $\U^n.$

 Let $D_i=\Fp(\Gamma T_i\cap \overline{\V}^n),$ $i=1,\ldots, k.$ These are closed in the metric topology of $\X,$ by \ref{remT}.
 
 By definition $\bigcup_{i=1}^k D_i\subs S,$ and
since $\Fp$ preserves the dimension $\dim_\R$  converting $\dim_\C$ into $2\cdot\dim_\R$  and taking into account that the missing components of $\Fp\inv(S)$ are inessential
we get
$$\dim_\R(S\setminus \bigcup_{i=1}^k D_i)\le \dim_\R S-2.$$
where the dimensions here are understood in the sense of \ref{loc-o-min}. Moreover,
recalling that the structure $(\U,\X,\Fp)$ is locally o-minimal we may deduce from the latter that in any small neighbourhood $V$ of a point of $S$
$\bigcup_{i=1}^k D_i\cap V$ contains an open subset of the analytic set $S\cap V.$ But  also $\bigcup_{i=1}^k D_i\cap V$ is closed in $V.$ It follows that $\bigcup_{i=1}^k D_i\cap V=S\cap V$ and hence
$\bigcup_{i=1}^k D_i=S.$

This immediately implies that every components of $\Fp\inv(S)$ intersecting $\V^n$ is essential. 

Note that $S$ is connected since it is an algebraically irreducible subvariety of a complex variety $\X^n.$
Moreover, it will stay connected if we remove from it a subset of real dimension $\dim_\R S-2.$ It follows
that  $\bigcup_{i=1}^k D_i$ is connected and moreover one can get from $D_1$ to any $D_i$ by a chain 
$D_1=D_{i_1},\ldots,D_{i_m}=D_i$ such that $\dim_\R(D_{i_j}\cap D_{i_{j+1}})\ge \dim_\R S-1,$ for
$1\le j<m.$ It follows that $\dim_\C (\gamma_jT_{i_j}\cap T_{i_{j+1}})=\dim_\C S$ for some gluing $\gamma_j\in \Gamma^n,$ for $j=1,\ldots,m-1.$   The statement of Lemma follows. $\Box$   
\epk
\bpk \label{A} {\bf Lemma.} {\em  In notation of \ref{setting} for any  components $T_1,T_2$ of
$\Fp\inv(S)$ there is a $\gamma_{12}\in \Gamma^n$ such that $T_2=\gamma_{12} T_1.$}

{\bf Proof.} We may assume that $T_1$ is intersecting $\V^n$ and $T_2$ intersecting $\alpha \V^n$ for some $\alpha\in \Gamma^n.$ Now $\alpha\inv T_2$ intersects $\V^n$ and using Lemma~\ref{B} we get a required $\gamma_{12}.$ $\Box$ 
\epk

We will need the following.
\bpk \label{fact} {\bf Fact.} (Special case of Theorem 12.5 of \cite{PS2}) {\em Let $G\subseteq \C^n$ be a semi-algebraic open set and $X\subseteq G$ an irreducible complex analytic subset which is also semi-algebraic. Then there is a complex algebraic subset $X^{\rm Zar}\subseteq \C^n$ such that $X$ is an irreducible analytic component of the set  $X^{\rm Zar}\cap G.$}

Note that by \ref{f01} there are finitely many irreducible analytic components of the set  $X^{\rm Zar}\cap G.$

\epk
\bpk \label{ZcapU}{\bf Lemma.} {\em Given a weakly special $S\subs \X^n,$ for each analytic component  $T\subs \Fp\inv(S)$ there is
 a  Zariski closed Zariski irreducible subset $Z\subseteq \C^{mn}$ such that 
\begin{itemize}
\item[(i)]   $\dim_\C S=\dim_\C Z;$
\item[(ii)] 
 $T$ is an irreducible analytic component of the set  $Z\cap U;$  
 \item[(iii)]
 $$\Fp\inv(S)\cap Z= \bigcup_{\gamma\in \mathrm{St}(Z)}\gamma\cdot T,$$ 
for   $\mathrm{St}(Z)=\{\gamma\in \Gamma^n: \gamma\cdot (Z\cap U)=Z\cap U\}.$
 
\end{itemize}
     } 
     
{\bf Proof.}     We may assume $T$ is essential in $\V^n.$ By \ref{R} $T\cap \V^n$ is semi-algebraic and thus, for some open semi-algebraic $G\subs \U^n,$ $\dim_\R T\cap G=\dim_\R T.$ By \ref{fact} there is a Zariski closed set $Z=T^{\rm Zar}$ such that $T\cap G=Z\cap G.$ Clearly, the minimal such $Z$ must be Zariski irreducible. Since the complex algebraic dimension coincides with the complex analytic dimension on Zariski closed sets, and the latter is local, we have 
              $$\dim_\C Z=\dim_\C Z\cap G=\dim_\C T\cap G=\dim_\C S.$$ 
 
 This gives us (i). This also implies (ii) since by irreducibility of $T$ we will have $T\subs Z\cap \U^n,$ and by equality of dimension  $T$ has to be a component in  $Z\cap \U^n.$  
 
 Now (iii) follows by \ref{A}. $\Box$
\epk
\bpk \label{lemma-def} Given a Zariski closed irreducible  $Z\subs \C^{mn}$ we can rearrange
the analytic irreducible decomposition of $Z\cap \U^n$ so that
\be \label{invdec}Z\cap \U^n=\bigcup_{i\in \N} R_i \ee
where $$R_i=\bigcup_{\gamma\in \mathrm{St}(Z)}\gamma\cdot T_i,$$
for some $T_i,$   analytic irreducible component of $Z\cap \U^n.$ 
Obviously the decomposition (\ref{invdec}) is unique, up to enumeration. 

Clearly, $\mathrm{St}(R_i)=\mathrm{St}(Z).$

We call the $R_i$'s {\bf invariant analytic components of $Z\cap \U^n.$} 
\epk

\bpk \label{wcsdef}
{\bf Definition.} Given a weakly special $S\subset \X^n$ we call an analytic subset $\check{S}\subs \U^n$
  {\bf weakly co-special set associated with $S$} if $\check{S}= \Fp\inv(S)\cap Z$ for some Zariski closed Zariski irreducible subset $Z\subseteq \C^{mn}$ such that 
  $\dim_\C  \Fp\inv(S)\cap Z=\dim_\C \Fp\inv(S)=\dim_\C Z.$

  By \ref{ZcapU}(iii) and \ref{lemma-def}  $\check{S}$ is an invariant analytic components of $Z\cap \U^n.$
  
  We call $Z$ the Zariski closure of $\check{S}.$

\epk
\bpk  \label{wcsp0}{\bf Lemma.} {\em Given  a weakly co-special $\check{S}\subs \U^n$ and its Zariski closure $Z$ there exists an analytic set $\hat{S}\subs  \U^n$ complementing $\check{S}$ to $Z\cap \U^n,$ that is 
$$\check{S}\cup \hat{S}=Z\cap \U^n,\ \ \ \dim_\C(\check{S}\cap \hat{S})<\dim_\C \check{S}.$$}

{\bf Proof.} Following \ref{ZcapU} take for $\hat{S}$ the union of all the analytic irreducible components of $Z\cap \U^n$ which are not subsets of $\check{S}.$ $\Box$
\epk 
\bpk \label{wcsp1}{\bf Lemma.} {\em For a weakly special $S\subset \X^n$ and an associated weakly co-special  $\check{S}$
we have the decomposition
$$ \Fp\inv(S)=\bigcup_{\gamma\in \Gamma^n}\gamma\cdot\check{S}. $$
The set of (distinct) components $\gamma\cdot \check{S}$ are in bijective correspondence with
the cosets $\Gamma^n/{\mathrm{St}(\check{S})}.$  
}

{\bf Proof.} Immediate from \ref{A}. $\Box$
\epk
\bpk
{\bf Examples.} 1. For $\X$ an Abelian variety of dimension $g$ the weakly co-special subsets of $\U^n=\C^{ng}$ are cosets of $\C$-linear subspaces $L\subseteq  \C^{ng}$ such that $L+\Lambda^n$ is closed in  $\C^{ng}.$

2. For a Shimura variety $\C$ (the modular curve) the weakly co-special subsets of the upper half-plane $\U$ are just points and the weakly co-special subsets of $\U^2$ are the graphs of the maps $x\mapsto gx,$ for $g\in \mathrm{GL}^+(\Q).$ 

3. Let $\X=\C^\times\setminus \{ a\}$ be the complex torus $\C^\times$ with a 
point $a$ removed.  
Let $\U=\C\setminus \{ \ln a +2\pi i k: k\in \Z\}$ and $\Fp= \exp$ restricted to $\U.$ This satisfies all the condition \ref{f00}.

The weakly co-special subsets $\check{S}$ of $\U^n$ are the intersection of $\Q$-linear subspaces of $\C^n$ with $\U^n.$ Unlike Example 1 above  we can have here $\pr \check{S}$ not constructible. For example, for
$$\check{S}=\{ \la x,y\ra\in \U^2: y=2x\}$$ one has, for $\pr:\la x,y\ra\mapsto x,$
$$ \pr \check{S}= \U\setminus  \{ \frac{1}{2}\ln a +\pi i k: k\in \Z\}.$$ 
\epk
\bpk \label{wcsint} {\bf Lemma.} {\em  Suppose $\check{S}_1$ and $\check{S}_2$ are weakly co-special. Then 
$$\check{S}_1\cap \check{S}_2
=\bigcup_{i=1}^m\check{R}_i$$
for some weakly co-special $\check{R}_1,\ldots,\check{R}_m.$}

{\bf Proof.} By definition     $\check{S}_1$ and $\check{S}_2$ are invariant analytic components of    $\Fp\inv(S_1)\cap Z_1$ and  $\Fp\inv(S_2)\cap Z_2,$ correspondingly, for some Zariski irreducible $Z_1$ and $Z_2.$

Let \be \label{ws1}S_1\cap S_2=\bigcup_{i=1}^k{P_i}\ee be the decomposition into weakly special subsets, see \ref{f1}.

We will have correspondingly
\be \label{ws2}\Fp\inv(S_1)\cap \Fp\inv(S_2)=\bigcup_{i=1}^k \Fp\inv(P_i).\ee
It follows that the irreducible components  of the $ \Fp\inv(P_i)$ are exactly the
irreducible components of  $\Fp\inv(S_1)\cap \Fp\inv(S_2).$

By definition $$\check{S}_1\cap \check{S}_2=\Fp\inv(S_1)\cap \Fp\inv(S_2)\cap Z_1\cap Z_2,$$
so using the decomposition $$Z_1\cap Z_2=\bigcup_{j=1}^mY_j$$
into Zariski irreducible components, we get
\be \label{ws2+} \check{S}_1\cap \check{S}_2=\Fp\inv(S_1)\cap \Fp\inv(S_2)\cap \bigcup_{j=1}^mY_j.\ee

 Combining (\ref{ws2+}) with (\ref{ws2})
 \be \label{ws3} \check{S}_1\cap \check{S}_2=\bigcup_{i=1}^k \Fp\inv(P_i)\cap \bigcup_{j=1}^mY_j= \bigcup_{i,j} \Fp\inv(P_i)\cap Y_j. \ee
 
 Now note that for each $j$ the analytic subset $\check{S}_1\cap \check{S}_2\cap Y_j$ of $Y_j\cap \U^n$ can be complemented (see \ref{wcsp0}) by the analytic subset 
 $(\check{S}_1\cap \hat{S}_2 \cup \hat{S_1}\cap\check{S_2}\cup \hat{S_1}\cap \hat{S_2})\cap Y_j$ to $Y_j\cap \U^n.$ This implies that any irreducible component $T$ of $\check{S}_1\cap \check{S}_2$ is also a component of $Y_j\cap \U^n.$ It follows that $\dim_\C T=\dim_\C Y_j.$ It further implies that
 $$\dim_\C  \Fp\inv(P_i)\cap Y_j=\dim_\C  \Fp\inv(P_i)= \dim_\C   Y_j.$$

 Denote $$\check{R}_{i,j}:= \bigcup_{i,j} \Fp\inv(P_i)\cap Y_j.$$
 These are weakly co-special subsets giving the required components of  $\check{S}_1\cap \check{S}_2.$ $\Box$
 
\epk
\bpk \label{wcspr} {\bf Lemma.} {\em  Given a weakly co-special  $\check{S}\subs \U^n,$ $0\le m<n$ and a projection $\pr: \U^n\to \U^{m},$ we have
$$\check{R}\supseteq \pr \check{S}\supseteq \check{R}\setminus \bigcup_{\gamma\in \Gamma^{m}}\bigcup_{i=1}^k\gamma\cdot\check{P_i},$$
for some weakly co-special $\check{R}$ and $ \check{P}_1,\ldots,  \check{P}_k\subs \check{R}\subs \U^{m},$ $\dim_\C \check{P}_i<\dim_\C \check{R}.$  }

{\bf Proof.} By definition $\Gamma^n\cdot\check{S}=\Fp\inv(S),$ $\Fp(\check{S})=S$ for some weakly special $S\subs \X^n.$ It follows from  \ref{f1}(ii) that $R\setminus \bigcup_{i=1}^k P_i\subs  \pr S\subs R$ for some weakly special $R$ and $P_1,\ldots, P_k\subsetneq R\subs \X^{m}.$ Let $\check{R}$ be a weakly co-special set associated with $R,$ that is $\Fp(\check{R})=R.$ 

Since $\pr(\Fp\inv(A))=\Fp\inv(\pr(A))$ for every $A\subs \X^n,$ we have 
$$\Gamma^{m}\cdot\check{R}\supseteq \pr(\Gamma^n\cdot \check{S})\supseteq \Gamma^{m}\cdot\check{R}\setminus \Fp\inv(P_1\cup\ldots\cup P_k).$$
It follows that $$\check{R}\supseteq \pr \check{S}\supseteq \check{R}\setminus \Fp\inv(P_1\cup\ldots\cup P_k).$$
Let $\check{P}_i$ be weakly co-special sets associated with the $P_i$'s. Then
by definition $\Fp\inv(P_i)=\Gamma^{m}\cdot \check{P}_i$ and we get the required.
$\Box$
\epk 
\bpk \label{wcsdm} {\bf Lemma.} {\em  Given a weakly co-special  $\check{S}\subs \U^n,$ $0\le m<n,$ a projection $\pr: \U^n\to \U^{}$  and a number $d=\dim_\C \check{S}(a),$ where $\check{S}(a)$ is a fibre over the point $a\in \pr \check{S}$
of the minimum dimension,   
define
$$\check{S}^{(d)}=\{ b\in \pr \check{S}: \dim_\C \check{S}(b)= d\}.$$ 
Then there are weakly co-special  $ \check{P}_1,\ldots,  \check{P}_k\subset \pr \check{S}$ each of dimension less than $\dim_C \pr \check{S}$ such that
$$\check{S}^{(d)}\supseteq \pr \check{S}\setminus \bigcup_{\gamma\in \Gamma^{m}}\bigcup_{i=1}^k\gamma\cdot\check{P_i}.$$
}

{\bf Proof.} As in the proof of \ref{wcspr} consider the  weakly special $S\subs \X^n,$
$S=\Fp(\check{S}).$
 By  \ref{f1}(iii)  ${S}^{(d)}\supseteq \pr S\setminus \bigcup_{i=1}^k P_i$ for some weakly special  $P_1,\ldots, P_k\subset \X^{m}$ of dimension less than that of $\pr S.$
 
Taking $\check{P}_i$ to be weakly co-special sets associated with the $P_i$'s, we get the required.
$\Box$
\epk

\bpk {\bf Definition.}  $\U_{S_w}$ is the structure with the universe $\U$ and the basic $n$-ary relations given by  the  weakly co-special subsets of $\U^n.$

We consider $\U_{S_w}$ a {\em topological structure} in the sense of \cite{Zbook},
with closed subsets defined as finite unions of  weakly co-special sets. 

By our definitions above, in particular \ref{dim}, there is a dimension notion $\dim_\C$ defined for all projective sets (i.e. the constructible sets and their projections).

\epk
\bpk We call a  weakly co-special subset $\check{S}\subs \U^n$ {\bf simple} if it is infinite and has no proper infinite  weakly co-special subsets.

In \cite{Zbook}, Ch.6 
an {\em analytic Zariski structure} has been defined. In the special case of irreducible one-dimensional
structure, which corresponds to a simple case here,
a combinatorial pregeometry has been defined and a closure operator $\cl$ introduced, see ibid. 6.3.  

\medskip

A more narrow but appropriate definition of an {\bf $\omega_1$-proper Zariski geometry} $\M$ is introduced by B.Elsner in \cite{Elsn}. This definition requires that $\M$ is $\omega_1$-compact (that is whenever
$\{ A_i:i\in I\}$ is a finitely consistent countable family of constructible sets, $\bigcap_{i\in I}A_i\neq \emptyset$) and the further
axioms are:
\begin{enumerate}
\item[(Z0)]   The topology on $\M$ and its Cartesian powers is Noetherian.
\item[(Z1)] For any closed $S\subs \M^{n+p}$ and a projection $\pr: \M^{n+p}\to \M^n$ there is a countable   
family $R_j,$ $j\in J,$ of closed subsets of $\M^n$ such that 
$\pr S\supseteq \overline{\pr S}\setminus \bigcup_{j\in J}R_j$ and $\dim R_j<\dim \pr S.$ 
\item[(Z2)] For any simple $Y\subs \M^p$ and  closed $S\subs \M^n\times Y$ there is a number $l$ such that
for every $a\in \M^n$ either the fibre $S(a)=Y$ or $|S(a)|\le l.$
\end{enumerate}

The last axiom (Z3) requires a {\em presmoothness} condition for simple subsets of $\M.$ We do not use it below.

\epk

Before we formulate and prove the  main theorem of this section we need the following result which may be of interest in its own right. The proof of it is obtained jointly with Y.Peterzil.
\bpk \label{kobiandI}{\bf Proposition.} {\em Let $a$ be a positive real number and $D\subset \C$  an open domain of the form $\{ z\in \C: -a< \mathrm{Re}(z)<a\ \& \ -a< \mathrm{Im}(z)<a\},$
$p:D\to \C$ a holomorphic injective function defined on  $D$ and in a small neighbourhood around $D,$ and suppose that
there is a semi-algebraic 4-ary relation $T(u_1,u_2,u_3,,u_4)$ on $\C$ such that
for any $x,y,z,w\in D$
$$w=x\cdot y+z\Leftrightarrow T(p(w),p(x),p(y),p(z)).$$
Then $p$ is a complex algebraic function on $D,$ that is there is a polynomial $f\in \C[z_1,z_2]$ 
such that the graph of $p$ is an irreducible analytic component of the analytic subset
$$\{ \la z_1,z_2\ra\in D\times \C: f(z_1,z_2)=0\} \subset D\times \C.$$}

{\bf Proof.} First we note that 
the structure $(D; <, p,\  w=x\cdot y+z)$ is definable in the canonical o-minimal 
structure $\R_{an}$  (where $<$  is defined on the interval $(a,b)$)

In the terminology of \cite{PS0} the ternary relation  $w=x\cdot y+z$ defines a {\em normal family of curves} on the interval $(-a,a)$ in  $(D; <,  w=x\cdot y+z).$ Hence 
by the main result of \cite{PS0} a field $(R_1,+,\cdot)$ definably isomorphic to the
 field $\R$ of reals is definable in $(D; <, w=x\cdot y+z).$ Since $(D; <, w=x\cdot y+z)$ is semi-algebraic,  $(D; <, w=x\cdot y+z)$ is bi-interpretable with  $(R_1,+,\cdot).$
 
 On the other hand, since the map $p: D\to C=:p(D)$ induces an isomorphism,
$$ (D; <, w=x\cdot y+z)\to (C; <, T(u_1,u_2,u_3,u_4))$$
(where $<$ on the right is defined on $p(-a,a)$)
a field $(R_2,+,\cdot)$ definably isomorphic to the
 field $\R$ of reals is bi-interpretable with $(C;  T(u_1,u_2,u_3,u_4)).$ Moreover, this bi-interpretation is given by the same formulas as in the first case and it extends the isomorphism $p$ between the two structures to an isomorphism
 $$p: (R_1,+,\cdot)\to (R_2,+,\cdot).$$
 
 Now we have two fields definable and definably isomorphic in $(D; <,p,\  w=x\cdot y+z).$ It is well-known (and easy to prove) that the only field-automorphism definable in an o-minimal expansion of $\R$ is the identity. Hence the isomorphism obtained by composing
 $i_1: \R\to R_1,$ $p: R_1\to R_2$ and $i_2\inv: R_2\to \R$ is the identity.   Since $i_1$ and $i_2$ are semi-algebraic interpretations, we get that $p:R_1\to R_2$ is semi-algebraic. But this is bi-interpretable with $p: D\to C.$ Hence $p$ is semi-algebraic.
 
 Finally, recalling that $p: D\to \C$   is holomorphic,  by \ref{fact} we get the required characterisation of $p.$ $\Box$  
  
\epk 

\bpk \label{fcor3} {\bf Theorem.} 

{\em 
(i) $\U_{S_w}$ is an  analytic Zariski structure.  More precisely, $\U_{S_w}$ 
is an $\omega_1$-proper Zariski geometry (not necessarily presmooth).

(ii) Let  $\check{\Y}\subs \U^n$ be a  simple weakly co-special set considered 
 a substructure of $\U_{S_w}.$
 Then  the corresponding weakly special set $\Y:= \Fp(\check{\Y})$ is simple and 
 the combinatorial geometry on $\check{\Y}$ is isomorphic to the combinatorial geometry on $\Y.$
 
(iii) $\check{\Y}$ satisfies the Trichotomy principle, i.e. the geometry on $\check{\Y}$  is either trivial, or linear, or of algebro-geometric type,   
in which case
 $\dim_\C Y=1$ and every Zariski closed subset of $Y^k,$ $k\in \N,$ is weakly co-special. Moreover, in this case 
   the restriction $\Fp$ on $\check{\Y},$   $\Fp: \check{\Y}\to \Y,$ is an algebraic map from every irreducible analytic component of $\check{\Y}$ onto the corresponding weakly special subset $\Y\subs \X^n.$ }

\smallskip

{\bf Proof.} As was noted above  $\U_{S_w}$ is a topological structure with a good dimension, in the sense of \cite{Zbook}.  The axioms on irreducible components and intersections follow from definition and \ref{wcsint}. The axiom on projections is given by \ref{wcspr}. The dimension of fibres condition is \ref{wcsdm}.
The rest follow from the fact that our formally analytic sets are actually complex analytic.


The topological structure $\U_{S_w}$ satisfies Elsner's axioms.

First note that $\U_{S_w}$ is $\omega_1$-compact. This is equivalent to the statement that if a constructible set $P$ is a union of countably many constructible subsets $P_i,$ $i\in \N,$ then 
$P=\bigcup_{i=1}^k P_i$ for some $k\in \N.$ The latter follows by induction on $\dim_\C P,$ using the fact that an irreducible analytic set  (the copies of which comprise a   weakly co-special subset) can not be represented as a countable union of closed subsets of smaller dimension.

(Z0) is given by \ref{wcsint} and (Z1) by \ref{wcspr}.  To see (Z2) one notes that the fibre $S(a)$ can be obtained by intersecting the weakly co-special sets $\U^n\times Y$ with $\{ a\}\times \U^p.$ Now we can use \ref{wcsint} which tells us that $S(a)$ is a finite union of weakly co-special subsets, in fact, checking the proof, the number of subsets is bounded by
$k\cdot m,$ where $k$ is the number of components in   
$(\X^n\times \Fp(Y))\cap  ( \{ \Fp(a)\}\times \X^p),$ and $m$ is the number of components in the intersection of Zariski closures  of $\U^n\times Y$ and $\{ a\}\times \U^p$ in
 $\C^{m(n+p)}.$ The first does not depend on $a$ and the second is just $\Fp(a)\times \C^{mp},$ hence $m$ does not depend on $a.$   But $k$ is also bounded from above as the number of components in a Zariski fibration. 
 
 It follows that either $S(a)$ is just $\{ a\}\times Y$ or  $S(a)$ is the union of at most
 $k\cdot m$ proper weakly co-special subsets of $\{ a\}\times Y$. Since $Y$ is simple the subsets are finite, so singletons. This proves (Z2) and finishes the proof of (i).
 
  

 (ii) We define {\em combinatorial dimension} $\mathrm{cdim}(A)$ for finite subsets $A\subset \check{\Y}$
following   \cite{Elsn}:
$$\mathrm{cdim}(a_1,\ldots,a_n)=\min \{ \dim \check{S}: \la a_1,\ldots,a_n\ra \in \check{S},\ \check{S}\subs_{\rm cl} \check{Y}^n\}.$$
The original definition in \cite{Zbook} and \cite{Elsn} assumes that $S$ runs among the projective subsets of $\U^n$ but as Elsner notes, under the given axioms (essentially (Z1))  we may assume that $S$ is closed. 

Combinatorial dimension gives rise to a closure operator and a pregeometry on $\check{Y}$
given by
$$b\in \cl(a_1,\ldots,a_n)\Leftrightarrow \mathrm{cdim}(a_1,\ldots,a_n)=\mathrm{cdim}(a_1,\ldots,a_n,b).$$

Since weakly special subsets $S$ of $\X^n$ and associated with them weakly co-special subsets $\check{S}$ of $\U^n$ are related by $\Fp(\check{S})=S$ we have
$$b\in \cl(a_1,\ldots,a_n)\Leftrightarrow \Fp(b)\in \cl(\Fp(a_1),\ldots,\Fp(a_n)),$$
where $\cl$ on the right is defined on $\Y$ by the corresponding condition.
This proves (ii)

(iii) The
  Trichotomy Theorem  follows from \ref{fcor2iii}.
 
It remains to prove that in the case the geometry is not locally modular,
 $\Fp: \check{\Y}\to \Y$ is an algebraic map. 
 
 Suppose $\Y$ is algebro-geometric, that is $\Y$ is an algebraic curve and every Zariski closed subset of $\Y^n$ is weakly special. We also have $\dim_\C\check{\Y}=1,$ since $\Fp$ preserves dimension.
 
Let $C\subset \check{Y}\cap \V^n$ be an open subset such that $\Fp:C\to \Fp(C)$ is bi-holomorphic.
  
 Let $f:\Y\to \C $ be given by a rational map. We may choose $C$ and  $f$ so that $f$ has no singularities on $\Fp(C)$ and so
 $f\inv$ exists as an algebraic function on $f(\Fp(C))\subset \C.$ We may assume $0\in  f(\Fp(C))$ and choose
  a neighbourhood $D\subset  f(\Fp(C))$ of $0$  of the form
 described in \ref{kobiandI}. We can even adjust $C$ so that $D=  f(\Fp(C)).$

 Now consider the holomorphic isomorphism $p=\Fp\inv\circ f\inv: D\to C.$  Since $\check{\Y}$ is an analytic subset of the complex manifold, we have an induced embedding $C\subset \C,$ so we may assume
 $p: D\to \C.$ Finally, we note that the image $T\subset \C^4$ under $p$ of the Zariski closed set 
 $$R=\{ \la w,x,y,z\ra\in D^4: w=xy+z\}$$    
is semi-algebraic, since the  Zariski closed subset $f\inv(R)$ of $Y^4$ is weakly special. 

We have now satisfied all the assumptions of \ref{kobiandI}. Hence $p$ on $D$ is algebraic. Hence $\Fp$ on
$C$ is algebraic. But then $\Fp$ is algebraic on the irreducible analytic component of $\check{Y}$ which contains $C.$  Since all  irreducible analytic components are conjugated by semi-algebraic transformations   
$\gamma$ (see \ref{lemma-def}) we proved (iii) and the theorem. $\Box$

\epk
\section{Special and co-special sets and points} \label{ss}

In this section we use the extra assumptions.

\medskip

{\bf Assumption D.} (i) the Zariski closure of $\U$ in $\C^m$ is defined over $\Q.$

(ii) $\X$ is defined over $\tilde \Q,$

\bpk {\bf Definitions.}
Call an  $S\subset \X^n$ {\bf strongly special}  if: 

(i) $S$  is weakly special, 

(ii) $S$ is defined over $\tilde \Q,$ 

(iii) for some 
  weakly co-special $\check{S}$ corresponding to $S$ the Zariski closure
 $\check{S}^\mathrm{Zar}$  is defined over $\Q.$

\medskip

We denote ${\cal S}_s$ the family of strongly special sets.


\medskip

Define ${\cal S},$ the family of {\bf special} sets,  to be the minimal family of sets containing ${\cal S}_s$ and closed under Cartesian products, 
intersections, Zariski closures of projections and taking irreducible components. A {\bf special point} is a singleton set which is special. 

Clearly, $${\cal S}_s\subseteq {\cal S}\subseteq {\cal S}_w.$$

We denote $\X_{S}$ the structure on $\X$ with $n$-ary relations given by special subsets of $\X^n.$ 

\epk


\bpk \label{ex2.6}{\bf Example.} 
Consider the algebraic torus $\C^*$ as $\X$ and  $\exp: \C\to \C^*$ as $\Fp.$ 

By Lindemann the trivial special point $1$ is the only strongly special point and $0$ a strongly co-special point. 

Weakly special sets will be cosets of tori. The strongly special sets are exactly tori, that is $0$-definable connected algebraic subgroups of $(\C^*)^n.$

Now the special sets by our definition are exactly Zariski closed  subsets $0$-definable in the multiplicative 
group $(\C^*;\cdot, 1).$ The set of points defined by the condition $x^k=1$ ($k$-torsion points) is definable by $\exists y\, x^k=y\ \& \ y=1,$ and so 
any torsion point is special. It is easy to conclude that special sets are exactly the torsion cosets of tori.

\epk
\bpk \label{fcor4} {\bf Proposition.} {\em $\X_{S}$ is a  Noetherian Zariski structure, if we consider any singleton in $\X$
to be closed.

(ii) The simple special sets  satisfy the Trichotomy.}


{\bf Proof.} Immediate from \ref{fcor2} and \ref{fcor2iii}.$\Box$ 
\epk

\bpk \label{Sh}{\bf Example.} Let $\X$ be a Shimura variety.
The main theorem of \cite{UY} characterises weakly special subvarieties in the sense of Shimura as weakly special in the sense above. 
 Shimura-special points are strongly special in our sense, and Shimura-special sets are those which are weakly special and contain a special point. It follows that special sets and special points in both senses are the same.   
\epk
\bpk {\bf Conjecture A.} Any weakly special set is definable with parameters in  $\X_{S}.$ Equivalently, every weakly special subset
$P\subseteq \X^n$
is a fibre of a special subset $S\subseteq \X^{n+m}$ under the  projection $\X^{n+m}\to \X^n.$
\epk

\bpk Note that the reference in \ref{Sh} confirms the conjecture A for Shimura varieties. 

\ref{ex2.6}  confirms the conjecture A in the case of the algebraic tori $\mathbb{G}_m(\C)^n.$

Special sets are also well understood for Abelian varieties and   the conjecture A can be confirmed in this class as well.

For the general mixed Shimura varieties to the best of our knowledge this is open but may be well within the reach of conventional methods.
\epk

\bpk
{\bf Conjecture B.} Given a special subset $S\subseteq \X^n,$ the set of special points in $S$ is Zariski dense in $S.$

\medskip

This important property of special points and sets in mixed Shimura varieties was pointed out to the author by E.Ullmo.

Using the fact that $\X_S$ has elimination of quantifiers it is easy to see that Conjectures A and B together are equivalent to the following model theoretic conjecture.

\medskip

{\bf Conjecture AB.} Let $\X_S^0$ be the substructure of $\X_S$ the universe of which is the set of special points. Then
 $$\X_S^0\preccurlyeq \X_S.$$

 \medskip
{\bf Conjecture Z-P.} {\em For any algebraic subvariety $V\subs \X^n$ there is a finite list of special subvarieties
$S_1,\ldots, S_m\subsetneq \X^n$
such that, given an arbitrary special subvariety $T\subset \X^n$ and an irreducible component $W$ of the intersection $V\cap T,$ either $\dim W=\dim V+\dim T-\dim \X^n$ (a typical case), or 
 $W\subs S_i$ for some $i=1,\ldots,m$ (in the atypical case $\dim W>\dim V+\dim T-\dim \X^n$).}

See comments in section~\ref{s1}, \ref{ZP}. 
\epk
\bpk {\bf Remark.} The validity of Conjecture AB implies that  $\X^0_S$ is a Zariski structure (without the extra requirement for singletons to be closed).  
\epk

\bpk {\bf Definition.} A {\bf co-special subset} of $\U^n$  is a weakly co-special subset $\check{S}$ corresponding to a special subset  $S\subseteq \X^n.$ 

Correspondingly, a {\bf period} is a point in $\Fp\inv(s)$ for $s\in \X$ a special point. 

We denote $\U_S$ the structure on $\U$ with $n$-ary relations given by co-special subsets of $\U^n.$ 
\epk
\bpk  {\bf Proposition.} {\em $\U_{S}$ is an $\omega_1$-proper Zariski structure, if we consider closed any singleton in $\U.$

}

{\bf Proof.} Immediate from \ref{fcor3}(i), since by definition the family of co-special sets are closed
under taking irreducible components of intersections. $\Box$
\epk

\bpk {\bf Definition.} Given $u_1,\ldots, u_n\in \U$ we define the {\bf special locus} of $\la u_1,\ldots, u_n\ra$
to be the minimal co-special subset $\check{S}$ of $\U^n$ containing $\la u_1,\ldots, u_n\ra$ and definable over $\tilde \Q$.
 Write in this case $\check{S}=\mbox{Sp.locus}(u_1,\ldots, u_n).$

We define a (combinatorial) dimension for tuples of points in $\U$ by setting, for $u_1,\ldots, u_n\in \U$
$$d_\mathrm{spec}( u_1,\ldots, u_n)=\dim_\C \mbox{Sp.locus}(u_1,\ldots, u_n).$$

\epk
\bpk {\bf Examples.} If $\Fp=\mathrm{id}$ and $\U=\X$ then special and co-special subsets of $\X^n$ and $\U^n$ are exactly the Zariski closed subsets defined over $\tilde \Q.$ In this case $d_\mathrm{spec}( u_1,\ldots, u_n)=\trd_\Q( u_1,\ldots, u_n).$

If $\Fp=\exp,$  $\U=\C$ and $\X=\C^*,$ then  co-special subsets of $\C^n$ are given by systems of equations of the form
$$m_1u_1+\ldots+m_nu_n=2\pi ik$$
for $m_1,\ldots,m_n,k\in \Z.$ These are defined over $\tilde \Q$ if and only if $k=0.$
Hence $$d_\mathrm{spec}( u_1,\ldots, u_n)=\ld_\Q(u_1,\ldots,u_n),$$
the $\Q$-linear dimension. 

\medskip

 Recall that Schanuel's conjecture can be stated as the following
 $$\trd_\Q(u_1,\ldots, u_n, e^{u_1},\ldots, e^{u_n})- \ld_\Q(u_1,\ldots,u_n)\ge 0.$$

  Now we have all ingredients to state the most general form of a Schanuel-type conjecture.

\epk
\bpk {\bf Conjecture D.}  Under the definitions above
$$\trd_\Q(u_1,\ldots, u_n, \Fp({u_1}),\ldots, \Fp({u_n}))- d_\mathrm{spec}(u_1,\ldots,u_n)\ge 0.$$

We do not know how exactly this conjecture is related to the Andr\'e - Grothendieck conjecture, which has a different form. However, all the instances of the latter conjecture discussed in \cite{An}, subsection 23.4, and in \cite{Be}  can be formulated as instances of Conjecture D.

\epk
\end{document}